\documentclass[11pt]{article}

\usepackage{amssymb,latexsym}
\usepackage{amsmath}
\usepackage{enumerate}

\title{On representations of star product algebras over cotangent spaces
       on Hermitian line bundles}

\author{{\bf
          Martin
          Bordemann$^a$\thanks{Martin.Bordemann@physik.uni-freiburg.de}~,
          \addtocounter{footnote}{1}
          Nikolai
          Neumaier$^a$\thanks{Nikolai.Neumaier@physik.uni-freiburg.de}~,
          Markus J.
          Pflaum$^b$\thanks{pflaum@mathematik.hu-berlin.de}~,}
          \\[3mm]
          {\bf
          Stefan
          Waldmann$^a$\thanks{Stefan.Waldmann@physik.uni-freiburg.de}
         } \\[1cm]
         {\begin {minipage} {7cm}
             \begin {center}
                 $^a$Fakult\"{a}t f\"{u}r Physik \\
                 Universit\"{a}t Freiburg \\
                 Hermann-Herder-Str. 3 \\
                 79104 Freiburg i.~Br., F.~R.~G
             \end {center}
         \end {minipage}
         \begin {minipage} {7cm}
             \begin {center}
                 $^b$Mathematisches Institut \\
                 Humboldt Universit\"{a}t zu Berlin \\
                 Unter den Linden 6 \\
                 10099 Berlin, F.~R.~G.
             \end {center}
         \end {minipage}
         }
         \vspace{1cm}
         }

\date{FR-THEP-98/24 \\[1mm]
      4.~December, 1998}

\newcommand{\LS} [1] {{#1} (\!(\lambda)\!)}
\newcommand{\CNP} [1] {{#1} \langle\!\langle \lambda \rangle\!\rangle}
\newcommand{\order} {\kappa}
\newcommand{\TinyW} {{\mbox{\rm \tiny W}}}

\newcommand{\starw} {\mathop{\star_{\mbox{\rm \tiny W}}}}
\newcommand{\wrep} {\varrho_{\TinyW}}
\newcommand{\Aw} {\mathcal A_{\mbox{\rm \tiny W}}}
\newcommand{\starwB} {\mathop{\star^B_{\mbox{\rm \tiny W}}}}
\newcommand{\wrepA} {\wrep^A}

\newcommand{\ad}{{\rm ad}}
\newcommand{\stars} {\mathop{\star_0}}
\newcommand{\srep} {\varrho_0}

\newcommand{\starsB} {\mathop{\star_0^B}}
\newcommand{\srepA} {\srep^A}
\newcommand{\rep} {\varrho}

\newcommand{\Ao} {{\mathcal A}_\order}
\newcommand{\staro} {\mathop{\star_\order}}
\newcommand{\repo} {{\rep_\order}}
\newcommand{\repAo} {{\rep_\order^A}}
\newcommand{\starBo} {\mathop{\star_\order^B}}
\newcommand{\repLo} {{\eta_\order^L}}
\newcommand{\repLAo} {{\eta_\order^A}}
\newcommand{\srepLA} {{\eta_{0}^A}}
\newcommand{\wrepLA} {{\eta_{\mbox{\rm \tiny W}}^A}}
\newcommand{\TAU} [3] {\tau^{{#1}}_{({#3},{#2})}}
\newcommand{\ormid} [2] {\gamma_\order ({#1},{#2})}
\newcommand{\Aj} {{\mathcal A}^{j}}
\newcommand{\half} {\frac{1}{2}}
\newcommand{\halfstars}{\bullet_0}
\newcommand{\halfstarw} {\bullet_{\mbox{\rm \tiny W}}}
\newcommand{\halfsrep} {\pi_0}
\newcommand{\halfwrep} {\pi_{\mbox{\rm\tiny W}}}
\newcommand{\starloc} {\star_{\mbox{\tiny \rm loc}}}
\newcommand{\Ds} {{\mathsf D}}
\newcommand{\Fdiff}[1]{{\mathsf F} \left( {#1} \right) }
\newcommand{\Exp}{{\rm Exp}}
\newcommand{\starf} {\star_{\mbox{\rm \tiny F}}}
\newcommand{\im} {{\mathrm i}}
\newcommand{\eu} {{\mathrm e}}

\newcommand{\Lie} {\mathcal L}
\newcommand{\cc} [1] {\overline {{#1}}}
\newcommand{\supp} {{\rm supp}}

\newcommand{\id} {{\sf id}}
\newcommand{\tr} {{\sf tr}}

\newcommand{\Hom} {{\mathsf {Hom}}}
\newcommand{\End} {{\mathsf {End}}}

\newcommand{\graph} {{\rm graph}}
\newcommand{\Op} {{\mathsf {Op}}}
\newcommand{\Opo} {{\mathsf {Op}_{\hbar,\order}}}
\newcommand{\Opw} {{\mathsf {Op}_{\hbar,\mbox{\rm \tiny W}}}}
\newcommand{\sym}{{\mathrm S}}
\newcommand{\csym} {{\mathrm S}_{\rm cpt}}
\newcommand{\Nopo} {N_\order^{\rm op}}
\newcommand{\brac}[1]{\left\langle{#1}\right\rangle}
\newcommand{\R}{{\mathbb R}}
\newcommand{\N}{{\mathbb N}}
\newcommand{\C}{{\mathbb C}}
\newcommand{\Z}{{\mathbb Z}}

\newcommand{\DL} {{{\mathsf {D}}^{L}}}
\newcommand{\degs} {{\mathrm {deg}_{s}}}

\renewcommand{\div}{{\sf div}}

\newenvironment{proof}{{\sc Proof:}}{{\hspace*{\fill} $\square$\\}}

\newtheorem{lemma} {Lemma} [section]
\newtheorem{proposition} [lemma] {Proposition}
\newtheorem{theorem} [lemma] {Theorem}
\newtheorem{corollary} [lemma] {Corollary}
\newtheorem{definition}[lemma] {Definition}
\newtheorem{example}[lemma] {Example}
\newtheorem{Remark}[lemma]{Remark}
\newenvironment{remark}{\begin{Remark}\rm}{\end{Remark}}
\numberwithin{equation}{section}

\begin{document}

\maketitle

\begin{abstract}
  For every  formal power series $B=B_0 + \lambda B_1 + O(\lambda^2)$
  of closed two-forms on a manifold $Q$ and every value of an ordering
  parameter $\kappa\in [0,1]$ we construct a concrete star product $\starBo$
  on the cotangent bundle $\pi : T^*Q\to Q$. The star product
  $\star_\order^B$ is associated to the formal symplectic form on $T^*Q$
  given by the sum of the canonical symplectic form $\omega$ and the
  pull-back of $B$ to $T^*Q$. Deligne's characteristic class of $\star_\order^B$
  is calculated and shown to coincide with the formal de Rham cohomology class
  of $\pi^*B$ divided by $\im\lambda$. Therefore, every star product on $T^*Q$
  corresponding to the Poisson bracket induced by the symplectic form
  $\omega + \pi^*B_0$ is equivalent to some $\starBo$. It turns out that
  every $\starBo$ is strongly closed. In this paper we also construct and
  classify explicitly formal representations of the deformed algebra as
  well as operator representations given by a certain global symbol
  calculus for pseudodifferential operators on $Q$. Moreover, we show that
  the latter operator representations induce the formal representations
  by a certain Taylor expansion. We thereby obtain a compact formula for
  the WKB expansion.
\end{abstract}

\newpage

\tableofcontents

\newpage

\section*{Introduction}
\label{IntroSec}
\addcontentsline{toc}{section}{\numberline{}Introduction}
In recent years the theory of deformation quantization has gained more and
more interest, not only because it provides a successfull approach to the
quantization problem (see e.g.~the work of {\sc Bayen} et al.~\cite{BFFLS78},
{\sc Fedosov} \cite{Fed:DQIT} and {\sc Kontsevich} \cite{Kon97b}),
but also because there are fascinating connections to analysis which
culminate in a deformation theoretical proof of the
Atiyah--Singer--Index Theorem (cf.~{\sc Fedosov}
\cite{Fed:DQIT} and {\sc Nest, Tsygan} \cite{NT95a,NT95b}). The main object
in the theory of deformation quantization is a formal deformation,
the so-called star product
$\star$, of the commutative algebra of all smooth complex-valued
functions ${\mathcal C}^\infty (M)$ on a given Poisson manifold
$(M, \{ \: , \: \})$,
such that in first order of the deformation parameter $\lambda$ the
$\star$-commutator equals a multiple of the Poisson bracket.
More precisely, $\star$ comprises a $\C[[\lambda]]$-bilinear associative
product
\begin{equation}\label{produitstar}
  \star: {\mathcal C}^\infty (M)[[\lambda]] \times
  {\mathcal C}^\infty (M)[[\lambda]] \rightarrow
  {\mathcal C}^\infty (M)[[\lambda]], \quad (f,g) \mapsto
  f \star g =  fg +\sum_{k=1}^\infty C_k (f,g) \lambda^k
\end{equation}
where the
$C_k: {\mathcal C}^\infty (M)\times {\mathcal C}^\infty (M)\rightarrow
{\mathcal C}^\infty (M)$ are assumed to be bidifferential operators
which satisfy $C_1(f,g)-C_1(g,f)=\im\{f,g\}$ and  $C_k(1,g)=C_k(f,1)=0$.
\vspace{2mm}

Although the existence and classification of star products on $M$ up to
equivalence is by now fairly well-understood,
there seem to be not very many results dealing with the problem of how to
define and construct some kind of operator representations for $\star$.
There are at least two principal choices.
The first, more algebraic point of view is to stay within the category of
modules over the ring $\C[[\lambda]]$ and regard representations as
$\C[[\lambda]]$-linear morphisms of the algebra
$({\mathcal C}^\infty (M)[[\lambda]],\star)$ to the algebra of
endomorphisms of some $\C[[\lambda]]$-module.
The second, more analytic approach consists in linearly mapping
${\mathcal C}^\infty (M)$ to operators on a suitable function space
such that the representation identity is satisfied `asymptotically',
a notion which depends on the topologies involved and which we shall make
more precise further down.

The latter, analytic notion of representation has already been described
in the work of {\sc Fedosov}, where a certain subspace
of ${\mathcal C}^\infty (M)[[\lambda]]$ is mapped to the space of
trace-class operators on a complex Hilbert space
(see e.g.~his book \cite[Chap.~7]{Fed:DQIT} for details).
{\sc Fedosov} also provides an existence proof for
his asymptotic operator representations which requires a certain integrality
condition for the symplectic form involved.

{\sc Bordemann} and {\sc Waldmann} have set up a {\it formal} representation
theory for deformation quantization
(see \cite{BW98a} and Appendix \ref{GNSApp} for a short
outline). There, a $\C[[\lambda]]$-analogue of the GNS construction
for $C^*$-algebras is introduced using $\C[[\lambda]]$-linear formally
positive functionals
${\mathcal C}^\infty (M)[[\lambda]]\rightarrow \C[[\lambda]]$ which are
shown to occur in a number of examples relevant to quantum theory.

The representation problem becomes more interesting in the class of cotangent
bundles $\pi:T^*Q\rightarrow Q$ equipped with the canonical 2-form
$\omega$.
Here, the space ${\mathcal C}^\infty (T^*Q)[[\lambda]]$ is directly related to
the space of symbols of pseudodifferential operators over the base manifold
$Q$, and calls for a comparison of formal representations and concrete
pseudodifferential operator realisations of star product algebras on $T^*Q$.
To this end, {\sc Pflaum} \cite{Pfl:NSRM,Pfl98c} has defined a family
$\star_\kappa$ of star products on $T^*Q$ parametrized by a so-called
ordering parameter $\kappa\in [0,1]$. This approach is based on a global
symbol calculus using the Levi-Civita connection $\nabla$ of a
fixed Riemannian metric on $Q$. The resulting star product $\star_0$,
also called the {\it standard star product} associated to $\nabla$,
arises as a particularly simple generalization of quantization by
standard ordering for polynomial symbols on $\R^{2n}$
(cf.~\cite[Eq.~(50.1)]{GS86}).
The products $\star_0$ and $\star_{1/2}$ have also been obtained
by {\sc Bordemann--Neumaier--Waldmann} in \cite{BNW97a,BNW97b}.
The starting point there is an arbitrary torsion-free connection
$\nabla$ on $Q$ from which $\star_0$ is defined by means of a
natural Fedosov--Weyl algebra construction. Moreover, $\star_{1/2}$ has been
obtained from $\star_0$ by an explicit equivalence transformation $N_{1/2}$
which depends on a fixed volume density $\mu$ on $Q$. Additionally,
formal GNS representations on the space
${\mathcal C}^\infty (Q)[[\lambda]]$ have been constructed.
In Section \ref{DefQuantSec} we will briefly recall the basic definitions
and some necessary precise statements of the above results.
\vspace{2mm}

In the present article we study again star products over
cotangent bundles $T^*Q$, but now equipped with a more general (formal) symplectic
form of the form $\omega_B = \omega + \pi^*B$,
where $B$ is a (formal) closed 2-form on the base manifold $Q$.
Thus we deal with a space of two-forms $\omega_B$ inducing all of the second
de Rham cohomology groups of $T^*Q$ and which
includes the Hamiltonian mechanics inside a magnetic field described
by $B$ (see e.g.~\cite [Chap.~3]{GS86}). Hereby, the integral of
$B$ over any closed two-dimensional surface in $Q$ can be interpreted
as the magnetic monopole charge inside that surface.

Our first main objective is to give a concrete coordinate-free construction
of star products on $T^*Q$ corresponding to the above symplectic forms,
and to classify their equivalence classes.

Secondly, we are interested in an explicit construction of
representations of the deformed algebras.
Due to the possible nonexactness of $B$ we can no longer expect to
represent them as formal or (pseudo-) differential operators on
${\mathcal C}^\infty (Q)[[\lambda]]$ resp.~${\mathcal C}^\infty (Q)$,
but rather on the space of all smooth sections
$\Gamma^\infty (L)$ of a complex line bundle $L$ over $Q$.
Let us explain this in more detail. A {\it formal representation} of
$({\mathcal C}^\infty (T^*Q)[[\lambda]],\star)$ on the line bundle $L$ then
is a $\C[[\lambda]]$-linear algebra homomorphism
   $$
      \varrho: \: {\mathcal C}^\infty (T^*Q)[[\lambda]] \rightarrow
      \End \big(\Gamma^\infty (L)\big)[[\lambda]],
      \quad f \mapsto \sum_{k=0}^\infty\,
      \lambda^k \varrho_k (f), \quad f \in {\mathcal C}^\infty (T^*Q) ,
   $$
where the $\varrho_k (f)$ are differential operators on $\Gamma^\infty (L)$.
Under an {\it operator representation} we understand a $\C$-linear
assignment of a pseudodifferential operator
$\Op_\hbar (a)$ on $\Gamma^\infty (L)$ for
every $\hbar \in \R^+ := \{ r\in \R|\: r > 0\}$ and every symbol
$a \in \sym^\infty (Q)\subset {\mathcal  C}^\infty (T^*Q)$ such that
for two symbols $a,b \in  \sym^\infty (Q)$ the asymptotic expansion
   $$
     \Op_\hbar (a) \cdot \Op_\hbar (b) \sim \Op_\hbar (ab) +
     \sum_{k=1}^\infty \, \hbar^k \Op_\hbar (C_k(a,b))
   $$
holds. This means (cf.~\cite{Hor:ALPDOIII}) that there exists a
decreasing sequence  $m_j \rightarrow -\infty$ such that
   $ \Op_\hbar (a) \cdot \Op_\hbar (b) - \big( \Op_\hbar (ab) +
   \sum_{k=1}^j \, \hbar^k \Op_\hbar (C_k(a,b)) \big)$
is a pseudodifferential operator of order $m_j$ for all $j \in \N$.
We now call the formal representation $\varrho$ {\it induced} by the
operator representation $\Op_\hbar$, if for every $m \in \N$
  $$
    \Op_{\hbar} (a) - \sum_{k=0}^m \, \hbar^k \varrho_k (a)
    = \hbar^{m+1} R_\hbar^m,
  $$
where $R_\hbar^m$ is a pseudodifferential operator which
depends continuously on $\hbar$ in the weak operator topology.  We will
denote  this situation by
$\Op_{\hbar} (a) = \sum_{k=0}^\infty \hbar^k \varrho_k (a) \mod \hbar^\infty$
and can then interpret  $\varrho(a)$ as a kind of Taylor expansion for
$\Op_\hbar (a)$. Note that this definition differs from the one given by
{\sc Fedosov}: the pseudodifferential operators we expect from quantum
mechanics such as Schr\"odinger operators will be in general unbounded and
not of the trace-class.
\vspace{2mm}

We have obtained the following main results.
Section \ref{QDiffSec} contains the heart of all later constructions.
It consists in a quantization or in other words formal abelian deformation of
the commutative group of all {\it fiber translating diffeomorphisms}. These
are diffeomorphisms of $T^*Q$ of the form $\zeta_x \mapsto \zeta_x + A_0(x)$,
where $A_0 \in \Gamma^\infty (T^*Q)$ is a one-form on $Q$ and
$\zeta_x \in T^*_xQ$. The outcome will be an isomorphism $\mathcal A$
interpolating between different  standard order
(or more generally $\order$-ordered)
star product algebras whose corresponding affine connections differ by
the one-form $A_0$.

After fixing a torsion-free affine connection $\nabla$ and a volume
density $\mu$ on $Q$ we provide in Section \ref{MagFieldSec} for every
value of the ordering parameter $\order \in [0,1]$ and  for every formal
power series of closed two-forms
$B= B_0 + \sum_{k=1}^\infty \, B_k \, \lambda^k$ a star product
$\starBo$ corresponding to the above formal symplectic form
$\omega_B = \omega + \pi^*B$. We call it
the $\order$-{\it ordered star product} associated to
$\nabla$, $\mu$ and $B$.
Hereby we first create the star product $\star_\order$ out of the standard
star product on $T^*Q$ by means of a global equivalence transformation
$N_\order$. Then we use the above-mentioned quantized fiber translations
along local potentials of $B$ to define local star products corresponding to
$\omega_B$.
It turns out that there are no obstructions to glue these local star
products together, so we obtain a global deformation $\starBo$ having the
desired properties. Moreover, applying the results of
{\sc Bertelsson--Cahen--Gutt} \cite{BCG97}, {\sc Deligne} \cite{Del95},
{\sc Fedosov} \cite{Fed:DQIT}, {\sc Gutt--Rawnsley} \cite{GR98} and
{\sc Weinstein-Xu} \cite{WX97} on the isomorphy classes of formal
deformations over symplectic manifolds, we show that every star product
on $(T^*Q, \omega +\pi^* B_0)$ is equivalent to some $\starBo$
the Deligne class of which is given by the de Rham class of $\pi^*B$
divided by $\im\lambda$. Another
nice property of the products $\starBo$ is that they are strongly closed,
i.e.~that the integral over the Liouville form provides a formal
{\em trace} for the algebra ${\mathcal C}^\infty (M)[[\lambda]]$.

In the following we construct {\it formal} representations for the star
products $\starBo$ explicitly. For exact $B=dA$ this can be done easily
by means of the global quantized fiber translating map along
$A$ (Section \ref{RepMagFieSec}).  Moreover, we arrive at a
compact closed formula for the WKB-expansion, and show afterwards
in Section \ref{RepHalDensSec} that $\star^B_0$ can be canonically
represented on the space of half-densities where the `magnetic field' $B$
is now a certain multiple of the trace of the curvature tensor.
In Section \ref{ABSec} we succeed in classifying the isomorphy classes of
representations associated to vector potentials of $B$
by de Rham cohomology classes modulo integral cohomology classes.
The underlying one-forms induce additional terms in the covariant
derivative on the trivial complex line bundle over $Q$ and thus may
create nontrivial holonomy.
In particular for vanishing $B$ this gives a possible interpretation of
Aharonov--Bohm like effects as known from physics in deformation quantization.

For nonexact $B=\lambda B_1 + O(\lambda^2)$ we can explicitly construct formal
representations of the deformed algebra on the space of all smooth sections
of a complex line bundle $L$ over $Q$ with Chern class  $[\frac{1}{2\pi} B_1]$
(Section \ref{MagMonoSec}). Hence, the two-form $B$ has to satisfy a certain
integrality condition which is also known in geometric quantization.

In the rest of our paper we consider concrete pseudodifferential operator
representations in the above defined sense.
A considerable part of Section \ref{SymbSec} is concerned about a global
pseudodifferential and symbol calculus on manifolds which
essentially goes back to the work of {\sc Widom} \cite{Wid:CSCPO,Wid78}
(cf.~also \cite{Saf:POLC,Pfl:NSRM}).
We introduce  a somewhat more general symbol calculus than there and
construct for every  $a \in \sym^\infty (Q)$ and every value of
the ordering parameter $\order \in [0,1]$  an operator
$\Op_{\hbar,\kappa} (a)$ on $\Gamma^\infty (L)$ which can be interpreted
as the $\order$-{\it ordered pseudodifferential operator} associated to $a$.
As the main result of Section \ref{SymbSec} we prove that the map
$\Op_{\hbar,\kappa}$ comprises an operator representation for the
star product $\star_\kappa^B$, and that the previously defined
formal representations are in fact induced by $\Op_{\hbar,\kappa}$.

Finally, the algebras ${\mathcal C}^\infty (T^*Q)$,
$\End_{\C[[\lambda]]} \big( \Gamma^\infty (L)[[\lambda]]\big)$ and
$\Psi^\infty (L)$ all carry a $^*$-structure depending on the choice of the
volume density $\mu$ and the Hermitian structure of $L$, which
is given respectively by complex conjugation, formal adjoint, and
(formal) operator adjoint. By spectral theoretical reasons and reasons
of applicability in quantum mechanics one is particularly interested
in star products and representations thereof which preserve the
$^*$-structures or in other words which fulfill
\begin{equation} \label{StarInvStructure}
  \cc f \star \cc g = \cc {g\star f} \quad \mbox{and} \quad \varrho(\cc f ) =
  \varrho (f)^*
\end{equation}
for all $f,g \in {\mathcal C}^\infty (T^*Q)$.
Now, this is the essential property of a star product of Weyl type
as it is defined precisely in Section \ref{DefQuantSec}.
Throughout our paper we have put much attention to this particularly
important case and proved that for the value $\order = 1/2$ of the ordering
parameter the star products $\star_\order^B$ and its representations
have indeed the property (\ref{StarInvStructure}).
Again thinking of $C^*$-algebras we also asked the question in Section
\ref{GNSSymbSec}, whether the representations $\varrho_{1/2}$ and
$\Op_{\hbar,1/2}$ are induced by a (formal) GNS construction.
Surprisingly, we succeeded not only to interpret the formal representations
of $\star_{1/2}^B$ as particular GNS representations, where each transversal
section in $\Gamma^\infty(L)$ serves as a cyclic vector, but could also
prove a kind of GNS theorem for $\Op_{\hbar,1/2}$.
\section*{Acknowledgements}
The authors are most indebted to S.~Goette for giving us the hint
to look for transversal sections. M.B.~and S.W.~would like to
thank the Mathematisches Institut of the Humboldt University of Berlin
and the SFB 288
for supporting a research stay in December 1997 out of which parts of this
paper have grown. M.P.~acknowledges support of the Fakult\"at f\"ur Physik
in Freiburg, where work on this paper has been finished.
\section{Preliminaries}
\label{NotationSec}
In this section, which is more of a technical nature and for the convenience
of the reader, we set up the notation and introduce some differential geometric
material needed in the sequel. \vspace{2mm}

Under $Q$ we will always understand a smooth $n$-dimensional manifold,
the configuration space, and under $\pi: T^*Q \to Q$ its cotangent bundle.
The canonical one-form on $T^*Q$ will be denoted by $\theta$,
the induced symplectic form by $\omega = - d \theta$, and
the Liouville form by $\Omega = \frac{1}{n!} (-1)^{[n/2]} \omega^n$.
The equation $i_\xi \omega = - \theta$ defines a unique vector field $\xi$ on
$T^*Q$ called the Liouville vector field.
We regard $Q$ as naturally embedded in $T^*Q$ as its zero section and
often write $\iota : Q \hookrightarrow T^*Q$ for that embedding.
Sometimes we will make use of local coordinates $q^1, \ldots, q^n$ on $Q$;
the induced canonical coordinates for $T^*Q$ will be denoted by
$q^1, \ldots, q^n, p_1, \ldots, p_n$.
If not stated otherwise functions and tensor fields on $Q$ are
assumed to be complex valued.

In this article $\mu$ will always be a volume density on $Q$, i.e.~a
smooth positive density $\mu \in \Gamma^\infty (|\!\bigwedge^n\!|\, T^*Q)$.
The induced densities on $T_xQ$ and $T^*_xQ$, $x \in Q$
are denoted by $\mu_x$ resp.~${\mu_x}^{\!\!\!*}$. Moreover, $\nabla$ always
means a torsion-free connection for $Q$.  By $E\rightarrow Q$ we denote a
complex vector bundle over $Q$ with Hermitian metric
$\brac{\cdot,\cdot}^E$, and by $\nabla^E$ a metric connection on $E$.

The volume density $\mu$ and the affine connection $\nabla$ determine
a one-form $\alpha_\mu \in \Gamma^\infty (T^*Q)$ by
\begin{equation} \label{alphaDef}
  \nabla_X \mu = \alpha_\mu (X) \, \mu,
  \qquad X \in \Gamma^\infty (TQ).
\end{equation}
Similarly, the divergence $\div_\mu (X)$ of a smooth vector field $X$
with respect to $\mu$ can be defined by the equation
\begin{equation}
  {\mathcal L}_X \mu = \div_\mu (X) \, \mu,
\end{equation}
where $\Lie_X$ is the Lie derivative with respect to $X$. One  checks
easily that
\begin{equation}
\nonumber
  \div_\mu (X) = \alpha_\mu (X) + \tr \left( Z \mapsto \nabla_Z X\right).
\end{equation}
This equation suggests to extend the divergence $\div_\mu $ to a
differential operator $\div_\mu^E$ on the space
$\Gamma^\infty ( \bigvee TQ \otimes E) $ by setting
\begin{equation}
\begin{split}
    \div_\mu^E & ( X_1 \vee \cdots \vee  X_k \otimes e)
    = \sum_{l=1}^k X_1 \vee \cdots
          \stackrel {l}{\widehat{\cdots}} \cdots \vee X_k \otimes
          \left( \div_\mu (X_l) s +  \nabla^E_{X_l} s \right) \\
    & + \sum_{{j,l=1 \atop j \ne l}}^k
            \left(\nabla_{X_j} X_l\right) \vee X_1 \vee \cdots
            \stackrel {l}{\widehat{\cdots}} \cdots
            \stackrel {j}{\widehat{\cdots}}
            \cdots \vee X_k  \otimes e, \qquad
     X_1,\cdots, X_k \in \Gamma^\infty (TQ), \: e \in \Gamma^\infty (E).
\end{split}
\end{equation}
Hereby $\stackrel{l}{\widehat{\cdots}}$ means to
leave out the $l$-th term. The for our purposes essential property of
$\div_\mu^E$ is stated in the following lemma.
\begin{lemma}
\label{LemAdjointNabla}
  Let $X$ be a smooth real vector field on $Q$.
  Then for any Hermitian vector bundle $E\rightarrow Q$ with
  metric connection $\nabla^E$ the (formal) adjoint of $\nabla_X^E$
  with respect to the scalar product
  $\brac{e, f}_{E,\mu} = \int_Q \brac{e, f}^E \mu$, where
  $e, f \in \Gamma^\infty_{\rm cpt} (E)$, is given by
  \begin{equation} \label{AdjointNabla}
    \left( \nabla_X^E\right)^* e =
    - \nabla_X^E e - \div_\mu (X) e = - \div_\mu^E ( X \otimes e ).
  \end{equation}
\end{lemma}
\begin{proof}
  By Stokes' Theorem we have
  \begin{equation}
  \nonumber
    \begin{split}
      0 & = \int_Q {\mathcal L}_X \left( \brac{e, f}^E \mu \right)
        = \int_Q \left( \brac{\nabla_X^E e, f}^E
          + \brac{e, \nabla_X^E f}^E
          + \brac{e, f}^E \, \div_\mu (X) \right) \mu,
    \end{split}
  \end{equation}
  hence the claim follows.
\end{proof}

A second differential operator on $\Gamma^\infty ( \bigvee TQ \otimes E) $
needed in the sequel
is given by the symmetrized covariant derivative
$\mathsf D: \Gamma^\infty (\bigvee^k TQ \otimes E) \to \Gamma^\infty
(\bigvee^{k+1} TQ \otimes E)$. It is defined by
\begin{equation} \label{SymCovDef}
  \Ds T ( X_1 ,..., X_{k+1} ) = \sum_{\sigma \in S_{k+1}}
  \left(\left( \nabla_{X_{\sigma(1)}} \otimes \id + \id \otimes
  \nabla_{X_{\sigma(1)}}^E \right) T\right)
  \left( X_{\sigma(2)}, \cdots, X_{\sigma(k+1)} \right),
\end{equation}
where $T$ is a smooth section of $\bigvee^k TQ \otimes E$, and the $X_j$ are
smooth vector fields on $Q$ (cf.~{\sc Widom}\cite{Wid:CSCPO}).

Consider now the pulled back bundle $\pi^*E \to T^*Q$ then we denote by
$\mathcal P^k (Q;E) \subset \Gamma^\infty (\pi^*E)$ the sections
which are \emph{polynomial of degree $k$ in the momenta}, i.e.~those
sections ${\mathsf p}: T^*Q \rightarrow \pi^*E$ having the property
that for every $x \in Q$ the map
$\mathsf p|_{T^*_xQ}: T^*_xQ \to E_x$ is a polynomial of
degree $k$ with values in $E_x$. Moreover, let
$\mathcal P (Q;E) = \bigoplus_{k=0}^\infty \mathcal P^k (Q;E)$.
For the trivial bundle $E= Q \times \mathbb C$ we simply write
$\mathcal P^k(Q)$ or $\mathcal P^k$ instead of $\mathcal P^k (Q;E) $
and note that ${\mathsf p} \in \mathcal P^k(Q) \subset C^\infty (T^*Q)$ if
and only if $\Lie_\xi {\mathsf p} = k {\mathsf p}$.

To every contravariant symmetric tensor field
$T \in \Gamma^\infty (\bigvee^k TQ)$
we assign a smooth function $J(T)\in {\cal P}^k$ by defining
$J(T) (\zeta_x) = \frac{1}{k!}
\brac{\zeta_x \otimes  \ldots \otimes \zeta_x, T_x} $
for all $\zeta_x \in T^*_xQ$, $x \in Q$.
Then $J : \Gamma^\infty (\bigvee TQ) \to {\cal P} (Q)$
turns out to be an isomorphism of  $\mathbb Z$-graded commutative algebras.

Next let us denote by
$\Delta_\mu^E : \Gamma^\infty (\pi^*E) \to \Gamma^\infty (\pi^*E)$
the generalized vector-valued Laplacian
\begin{equation} \label {VectorLaplaceDef}
    \Delta_\mu^E =
    \sum_j \nabla_{(dq^j)^{\rm v}}^\pi
    \nabla_{(\partial_{q^j})^{\rm h}}^\pi +
    \nabla_{(\alpha_\mu)^{\rm v}}^\pi,
\end{equation}
where $(\cdot)^{\rm v} $ resp.~$(\cdot)^{\rm h}$ means the vertical
resp.~horizontal lift to $T^*Q$ induced by $\nabla$,
and $\nabla^\pi$ is the pull-back of $\nabla^E$ to $\pi^*E$ via $\pi$.
Note that the definition (\ref {VectorLaplaceDef}) is independent
of the used chart.
\begin {lemma}
For ${\mathsf p} \in \mathcal P (Q; E)$ one has
$\Delta^E_\mu {\mathsf p} = J \div^E_\mu J^{-1} {\mathsf p}$.
\end {lemma}

Let us now define a few important differential
operators on the cotangent bundle $T^*Q$.
First assign to any covariant symmetric tensor field
$\gamma \in \Gamma^\infty (\bigvee^k T^*Q)$ a fiberwise acting differential
operator $\Fdiff{\gamma}$ on $T^*Q$  of order $k \in \mathbb N$ by defining
\begin{equation}
\label{PgammaDef}
   \left(\Fdiff{\gamma_1 \vee \cdots \vee \gamma_k}f\right) (\zeta_x)
   = \left. \frac{\partial^k}{\partial t_1 \cdots \partial t_k}
   \right|_{t_1 = \ldots =t_k =0} f \left(
   \zeta_x  + t_1 \gamma_1 (x) + \ldots +  t_k \gamma_k(x) \right).
\end{equation}
Hereby $\zeta_x \in T^*_xQ$,
$\gamma_1,\ldots ,\gamma_k \in \Gamma^\infty (T^*Q)$
and $f \in {\mathcal C}^\infty (T^*Q)$.
Clearly (\ref {PgammaDef}) extends linearly to an injective
algebra morphism from $\Gamma^\infty (\bigvee T^*Q)$ into the algebra
of differential operators of ${\mathcal C}^\infty (T^*Q)$.
Note that $\Fdiff{\gamma}$ commutes with
$\Fdiff{\gamma'}$ for any
$\gamma, \gamma' \in \Gamma^\infty (\bigvee T^*Q)$.

In this article we will also make use of the Laplacian of a metric $G$ on
$T^*Q$ which is naturally induced by the a priori chosen torsion-free
covariant connection $\nabla$ on $Q$
(cf.~\cite{YanIsh:TCB} or \cite [App.~A]{BNW97a}).
The connection $\nabla$ induces a splitting of $T (T^*Q)$ into
its horizontal and vertical subbundle.
Now, if $X,Y \in \Gamma^\infty (TQ)$ and
$\alpha, \beta \in \Gamma^\infty (T^*Q)$ are smooth sections, we denote
by $X^{\rm h}, Y^{\rm h}$ the horizontal lifts of $X,Y$ and by
$\alpha^{\rm v}, \beta^{\rm v}$ the
vertical lifts of $\alpha, \beta$ to vector fields on $T^*Q$.
Then $G$ is defined by
$G (X^{\rm h}, Y^{\rm h}) := 0 =: G (\alpha^{\rm v}, \beta^{\rm v})$ and
$G (X^{\rm h}, \alpha^{\rm v}) := 2\alpha (X) =: G (\alpha^{\rm v},
X^{\rm h})$.
Obviously $G$ is an indefinite metric on $T^*Q$, so its Laplacian
(`d'Alembertian') $\Delta$ is well-defined.
One checks easily that in local canonical coordinates $\Delta$ has the form
\begin{equation}
\label{DeltaDef}
    \Delta = \sum_k
             \frac{\partial^2}{\partial q^k \partial p_k}
             + \sum_{j, k, l}
             p_l \left(\pi^* \Gamma^l_{jk}\right)
             \frac{\partial^2}{\partial p_j \partial p_k}
             + \sum_{j, k}
             \left(\pi^* \Gamma^j_{jk}\right)
             \frac{\partial}{\partial p_k} ,
\end{equation}
where $\Gamma^l_{jk}$ are the Christoffel symbols of $\nabla$.
Consequently, the relation
\begin{equation} \label{EquivLap}
  \Delta_\mu^Q = \Delta + \Fdiff{\alpha_\mu}
\end{equation}
holds, where $\Delta_\mu^Q$ is the generalized
vector-valued Laplacian $\Delta_\mu^E$ associated to the trivial line bundle
$E = Q\times \C \rightarrow Q$.

\section{Deformation quantization}
\label{DefQuantSec}
Let us now recall some essential notions and results from the
theory of formal deformation quantization.

A star product $\star$ on a symplectic manifold $M$ is called of
\emph{Vey type}, if the order of every bidifferential operator $C_k$
in each of its variables is lower or
equal to $k$. It is called of \emph{Weyl type}, if
$C_k (f,g) = (-1)^k C_k (g,f)$ and if in addition
$C_k$ is real (resp.~imaginary) for $k$ even (resp.~odd). A star
product of Weyl type statisfies in particular
$\cc {f \star g} = \cc g \star \cc f$, i.e.~complex conjugation is an
algebra anti-automorphism. Hereby we have formally set
$\cc \lambda = \lambda$. Two star products $\star$ and $\star'$
are called \emph{equivalent} if there exists a formal power series
$S = \id + \sum_{k=1}^\infty \lambda^k S_k$ of differential
operators $S_k: \mathcal C^\infty (M) \to \mathcal C^\infty (M)$
such that $S1 = 1$ and $S(f\star g) = S(f) \star' S(g)$ for all
$f, g \in \mathcal C^\infty (M)[[\lambda]]$.
\vspace{2mm}

In case $M$ is a cotangent bundle $T^*Q$ we impose some more structure on
a reasonable star product, that means on a star product which should be
relevant for quantum mechanics.

First one has a second type of `ordering' for star products on $T^*Q$
besides the Weyl type: one calls a star product $\star$ of
\emph{standard order} or \emph{normal order type}, if
$(\pi^* u) \star f = (\pi^* u) f$ for all
$u \in {\mathcal C}^\infty (Q)$ and $f \in {\mathcal C}^\infty (T^*Q)$.

Next let us consider the homogeneity operator
${\mathsf H} = \lambda \frac{\partial}{\partial \lambda} + \Lie_\xi$
on ${\mathcal C}^\infty (T^*Q) [[\lambda]]$.
A star product on $T^*Q$  is said to be {\it homogeneous}
(cf.~\cite {DL83a}), if $\mathsf H$ is a derivation.
This implies in particular the simple but important fact that
${\mathcal P}[\lambda]$ is a $\mathbb C[\lambda]$-algebra with respect to
a homogeneous star product (cf.~\cite [Prop.~3.7]{BNW97a}).
Hence it is possible to substitute the formal parameter $\lambda$ by a
real number $\hbar$. Speaking  in analytical terms  this
means that in the homogeneous case the convergence problem for a
star product can be trivially solved  for functions polynomial in
the momenta. Moreover, it has been shown in \cite {BNW97b} that
all homogeneous star products are strongly closed.

In \cite {BNW97a} a homogeneous star product of standard order type
$\stars$ was constructed by a modified Fedosov procedure using the
connection $\nabla$ on $Q$. Additionally, a  canonical representation
$\srep$ of the deformed algebra by formal power series of
differential operators on the space ${\mathcal C}^\infty (Q)[[\lambda]]$
of formal wave functions was constructed (cf.~\cite [Thm.~9] {BNW97a}).
Explicitly this representation is given by
\begin{equation}
\begin{split}
\label {SRepDef}
    \srep (f) u & =
          \iota^* (f \stars \pi^* u) =
          \sum_{l=0}^\infty \, \sum_{j_1,\cdots,j_l} \,
          \frac{(-\im\lambda)^l}{l!}
          \iota^*\left( \frac{\partial^l f}
          {\partial p_{j_1} \cdots \partial p_{j_l}} \right)
          \, \frac{1}{l!} \brac{\partial_{q^{j_1}} \otimes \cdots
          \otimes \partial_{q^{j_l}} , \Ds^l u} \\
          & = \sum_{l=0}^\infty \,
          \frac{(-\im\lambda)^l}{l!}
          \iota^*\left(  \Fdiff{ \Ds^l u } \, (f)  \right)
          = \iota^* \Fdiff{\exp(-\im\lambda\Ds)u}f, \quad
          f \in {\mathcal C}^\infty (T^*Q)[[\lambda]], \:
          u \in {\mathcal C}^\infty (Q)[[\lambda]].
\end{split}
\end{equation}

After fixing a volume density
$\mu \in \Gamma^\infty (|\!\bigwedge^n \!|\, T^*Q)$ on $Q$ one can pass
for every value of an \emph{ordering parameter} $\order \in [0,1]$
to a new star product $\staro$ on ${\mathcal C}^\infty (T^*Q)[[\lambda]]$
in the following way. Consider the one-form
$\alpha_\mu$ and let $N_\order = N_\order (\alpha_\mu) $ be the operator
\begin{equation} \label {NeumaierDef}
    N_\order (\alpha_\mu) =
    \exp \left( - \im \lambda \order \left(\Delta + \Fdiff{\alpha_\mu}
    \right) \right) = \exp \left( - \im \lambda \order  \Delta_\mu^Q \right) .
\end{equation}
Then one can define a scale of star products $\staro$ each equivalent
to $\stars$ by
\begin{equation}
\label{WeylDef}
    f \staro g = N_\order^{-1}
    \left(N_\order (f) \stars N_\order (g)\right).
\end{equation}
It turns out that all $\staro$ are homogeneous since $N_\order$
commutes with $\mathsf H$. Moreover, a representation $\repo$ of $\staro$
on ${\mathcal C}^\infty (Q)[[\lambda]]$ is given by
$\repo (f) = \srep (N_\order f)$. For the case $\order = 1/2$
one obtains the so-called \emph{Weyl star product}
$\starw = \star_{1/2}$ with corresponding
\emph{Schr{\"o}dinger representation} $\wrep = \varrho_{1/2}$
having the following two nice properties:
\begin{equation}
\begin{align}
\label {CCAnti}
    \cc{f \starw g} & = \cc g \starw \cc f , \\
\label {WrepAdj}
    \int_Q \cc u \: \wrep (f) v \;\mu & = \int_Q \cc{\wrep (\cc f) u}
    \: v \; \mu, \qquad u,v \in {\mathcal C}^\infty_{\rm cpt}
    (Q)[[\lambda]].
\end{align}
\end{equation}
In other words $\wrep$ is a $^*$-representation of the (formal) $^*$-algebra
$ ({\mathcal C}^\infty (Q)[[\lambda]],\starw)$ with complex conjugation as
its $^*$-involution. In particular $\starw$ is of Weyl type.
Finally, the representation $\wrep$ coincides with the formal GNS
representation (see Appendix \ref {GNSApp}) with respect to the
functional
$f \mapsto \omega_\mu (f) = \int_Q \iota^*f \; \mu$ defined on the
two-sided ideal ${\mathcal C}^\infty_Q (T^*Q) [[\lambda]]$,
where ${\mathcal C}^\infty_Q (T^*Q) =
\left\{ f \in {\mathcal C}^\infty (T^*Q) \; | \; \supp (\iota^* f)
\mbox{ is compact} \right\}$ (cf. \cite[Prop.~4.2]{BNW97b}).
\section{Quantization of fiber translating diffeomorphisms}
\label{QDiffSec}
In this section we will provide a method for
quantizing fiber translating diffeomorphisms.
To this end we first compute explicitly the $\staro$-product of
a homogeneous function of degree $0$ on $T^*Q$ with an arbitrary function.
\begin {lemma}
Let $\gamma, \gamma' \in \Gamma^\infty (\bigvee TQ)$. Then
\begin{eqnarray}
\label{NPCom}
    \big[\Delta + \Fdiff{\gamma'}, \Fdiff{\gamma}\,\big]
    & = & \Fdiff{\Ds\gamma},  \\
    \label {NkonjF}
    N_\order^{-1} \circ \Fdiff{\gamma} \circ N_\order
    & = & \Fdiff{ \exp (\im \order \lambda \Ds)\, \gamma } .
\end{eqnarray}
\end{lemma}
\begin{proof}
For example by computation in local coordinate it is easy to see that
(\ref {NPCom}) is true, if $\gamma$ has symmetric degree $0$ or $1$.
Then (\ref {NPCom}) follows for arbitrary symmetric degree, since $\mathsf F$
is an algebra morphism and $\Ds$ a derivation. Eq.~(\ref {NkonjF}) is
an immediate consequence of (\ref {NPCom}).
\end{proof}
\begin{proposition}
\label{Propad}
Let $f \in \mathcal C^\infty (T^*Q)[[\lambda]]$ and
$u \in \mathcal C^\infty (Q)[[\lambda]]$. Then for all $\order \in [0,1]$
the equalities
\begin{eqnarray}
\label{ustarf}
    \pi^*u \staro f
    & = & \Fdiff {\exp (\im \order\lambda \Ds) u} f, \\
    \label {fstaru}
    f \staro \pi^*u
    & = & \Fdiff {\exp (-\im (1-\order)\lambda \Ds) u} f
\end{eqnarray}
hold and consequently also
\begin {equation} \label{adstaro}
    \ad_\order (\pi^*u) = N_\order^{-1} \ad_0 (\pi^*u) N_\order
    = \Fdiff{\frac{\exp (\im\order\lambda\Ds)
                          - \exp (-\im(1-\order)\lambda\Ds)}{\Ds}
                          du} .
\end {equation}
\end {proposition}
\begin {proof}
Due to the equivalence of $\staro$ and $\stars$ via $N_\order$ and the
fact that $N_\order \circ \pi^* = \pi^*$, it suffices by (\ref {NkonjF})
to proof (\ref {ustarf}) and (\ref {fstaru}) for the case $\order = 0$.
Assuming $\order = 0$ the relation
(\ref {ustarf}) is trivially fulfilled since $\stars$ is of standard
order type. To prove (\ref {fstaru}) we can restrict to
$f, \pi^*u \in \mathcal C^\infty (T^*Q)$ without $\lambda$-powers.
Since $\stars$ is homogeneous we have in order $\lambda^l$
at least $l$ vertical derivatives acting all on $f$ since $\pi^*u$ is
homogeneous of degree $0$. On the other hand $\stars$ is of Vey type
(see \cite [Lem.~5] {BNW97a}) whence there are at most and thus
\emph{exactly} $l$ derivatives in fiber direction acting on $f$.
Moreover, the coefficient functions of the bidifferential operator in order
$\lambda^l$ can thus only depend on $Q$ due to the homogeneity.
Thus they are determined by their values at the zero section
which are known from (\ref {SRepDef}) proving the proposition.
\end {proof}
Note that the products (\ref{ustarf}) and (\ref{fstaru}) do not depend on
$\alpha_\mu$ and thus not on the choice of the volume density $\mu$.
Moreover, the above expressions for the commutators are well-defined
formal power series in $\lambda$ with differential operators as
coefficients. One observes that the star product $\star_1$
is of \emph{anti-standard order type} in an obvious sense.

Next one makes the simple but crucial observation that $\ad_\order (\pi^*u)$
only depends on the one-form $du$. Thus Eq.~(\ref {adstaro}) still makes sense
as formal series of differential operators, if $du$ is replaced by an
arbitrary formal one-form
$A = A_0 + \lambda A_1 + O(\lambda^2) \in \Gamma^\infty (T^*Q)[[\lambda]]$.
This motivates to define an operator
$\delta_\order [A]: {\mathcal C}^\infty (T^*Q) [[\lambda]] \to
{\mathcal C}^\infty (T^*Q) [[\lambda]]$ by
\begin{equation}
\label {deltaDef}
   \delta_\order [A] =
   \Fdiff{\frac{\exp ( \im \order \lambda \Ds) -
   \exp (- \im (1-\order) \lambda \Ds)}{\Ds} \; A} .
\end{equation}
Now, the following result is immediate.
\begin{proposition}
 The operator $\delta_\order [A]$ is a
 derivation of $\staro$ if and only if $dA = 0$, and an inner
 derivation if and only if $A$ is exact, that means if and only if
 $A = du$ with $u \in {\mathcal C}^\infty (Q)[[\lambda]]$.
 Additionally, the relation
 \begin{equation} \label {deltaWdeltaS}
    \delta_\order [A] = N_\order^{-1} \circ \delta_0 [A] \circ N_\order
 \end{equation}
 is true. Finally, one has in lowest order an expansion of the form
 $\delta_\order [A] = \im \lambda \Fdiff{ A_0} + O(\lambda^2)$.
\end{proposition}
In deformation quantization the Heisenberg equation of motion for a given
Hamiltonian $H$ induces a one-parameter group of star product automorphisms
provided the corresponding classical Hamiltonian vector field of $H$ has a
complete flow (see e.g.~\cite [App.~B] {BNW97b}). Although the operators
$\delta_\order [A]$ are  derivations or in other words infinitesimals of star
product automorphisms only in the case $dA = 0$, we consider time evolution
equations induced by these maps, i.e.~differential equations of the form
\begin{equation} \label {TimeEvolution}
    \frac{d}{dt} F(t) =
    \frac{\im}{\lambda} \delta_\order [A] \, F (t), \qquad
    F(0) = f,
\end{equation}
where $f \in {\mathcal C}^\infty (T^*Q)[[\lambda]]$ is an initial value.
First of all note that the form $A_0$ in the expansion
$A = A_0 + \lambda A_1 + O(\lambda^2)$ induces a one-parameter group of
diffeomorphisms $\phi_t : T^*Q \to T^*Q$  by
\begin{equation} \label {phitDef}
    \phi_t (\zeta_x) = \zeta_x - tA_0 (x).
\end{equation}
Hereby we assume $A_0$ to be real. Then $\phi_t$ is the
flow of the vector field determined by the differential operator
$- \Fdiff{A_0}$. Note that $\phi_t$ is symplectic if and only if $dA_0 = 0$.
Moreover, $\phi_t$ is a fiber translation, and
\begin{equation}
    \phi^*_t \circ \Fdiff{\gamma} = \Fdiff{\gamma} \circ
    \phi^*_t
\end{equation}
holds for all $\gamma \in \Gamma^\infty(\bigvee T^*Q)[[\lambda]]$.
Using this fact we obtain the following result.
\begin{theorem}
\label{EvolutionTheo}
Let $A \in \Gamma^\infty (T^*Q)[[\lambda]]$ with real $A_0$.
Then for every value of the ordering parameter $\order$ and every
initial value $f \in {\mathcal C}^\infty (T^*Q)[[\lambda]]$
Eq.~(\ref {TimeEvolution}) has a unique solution
$F : \R \rightarrow {\mathcal C}^\infty (T^*Q)[[\lambda]]$. It is given by
$F (t) = \Ao (t) f$, where the operator
${\mathcal A}(t) = \Ao (t) : {\mathcal C}^\infty (T^*Q)[[\lambda]]
\rightarrow {\mathcal C}^\infty (T^*Q)[[\lambda]]$ is defined by
\begin{equation}
\label{DeveloperS}
   \Ao (t) = \phi^*_t \circ  \exp \left(-t\Fdiff{
   \frac{\exp(\im\order\lambda\Ds) - \exp (- \im(1-\order)\lambda\Ds)}
   {\im \lambda \Ds}(A) - A_0 }\right),
\end{equation}
and has the following properties.
\begin{enumerate}
\item $\Ao (t)$ commutes with $\Fdiff{\gamma}$ and the
      pull-back $\psi^*$ for all $\gamma \in \Gamma^\infty (\bigvee TQ)$
      and all fiber translating diffeomorphisms $\psi$.
\item For all times $t \in \R$ one has
      \begin{equation} \label {DeveloperW}
         \Ao (t) =
         N_\order^{-1} \circ \mathcal A_0 (t) \circ N_\order .
      \end{equation}
\item If $A' \in \Gamma^\infty (T^*Q)[[\lambda]]$ is another one-form
      fulfilling $A'_0 = \cc{A'_0}$,  then the equality
      \begin{equation} \label {AswAsw}
         {\mathcal A} (t) \circ {\mathcal A}' (t') =
         \tilde{\mathcal A} (1) ={\mathcal A}' (t') \circ {\mathcal A} (t)
      \end{equation}
      holds for all $t, t' \in \mathbb R$, where
      ${\mathcal A}'$ resp.~$\tilde{\mathcal A}$ denotes the
      evolution operator for the one-form $A'$ resp.~$tA + t'A'$.
\item The map $t \mapsto \Ao (t)$ comprises a one-parameter group of
      $\staro$-automorphisms if and only if $dA = 0$.
\item If $A$ has the form
      $A = \lambda A_1$, then ${\mathcal A}(t)$ commutes with $\mathsf H$.
      In case $A = \cc A$,  then
      $\Aw (t) := {\mathcal A}_{1/2} (t)$ commutes with
      complex conjugation.
\end{enumerate}
\end{theorem}
\begin {proof}
It is a well-known fact that the Heisenberg equation (\ref {TimeEvolution})
has a unique solution for any initial value if the flow of the vector field
corresponding to the zeroth order is complete (see e.g.~\cite {BNW97b}).
Thus it remains to show that (\ref {DeveloperS}) is indeed a solution. But
this is easily obtained by differentiation. From (\ref {DeveloperS}) the
first part follows directly as well as part ii.) which is no surprise due
to (\ref {deltaWdeltaS}). Part iii.) is computed using the above stated
properties of fiber translating diffeomorphisms as well as the fact that
all involved operators commute. Part iv.) is known in the general case
(see e.g.~\cite {BNW97b}) but can also be verified explicitly using
(\ref {DeveloperS}). Moreover, if $A = \lambda A_1$ then
$\frac{\im}{\lambda}\delta_\order[A]$ commutes with $\mathsf H$ and so does
$\Ao (t)$. This can be derived from (\ref {TimeEvolution}) or explicitly from
(\ref {DeveloperS}). Additionally, (\ref {DeveloperS}) implies for
$A = \cc A$ and $\order = 1/2$ that the operator
$\Aw (t) = \mathcal A_{1/2} (t)$ commutes with complex conjugation.
\end {proof}
Note that in general ${\mathcal A}$ is not homogeneous,
i.e.~that it does not commute with ${\mathsf H}$.
But nevertheless one has the relation ${\mathcal A} \circ \pi^* = \pi^* $
for arbitrary $A$, see also Lemma \ref {DeveloperHComLem}.

The above formulas suggest that we call $\Ao (t)$ the
\emph{quantization} of the fiber translating diffeomorphism
$\phi_t$ for the value $\order$ of the ordering parameter.
\begin{example}
\rm
Let $X \in \Gamma^\infty(TQ)$ be a vector field. Then $J(X)$
is a function on $T^*Q$ linear in momentum.
Hence one computes for $A \in \Gamma^\infty(T^*Q)[[\lambda]]$ that
$\delta_\order [A] \, J(X) = \im\lambda \pi^* (A(X))$ and
\begin{equation} \label {MiniCouple}
    \Ao (t) J(X) = J(X) - t \pi^*(A(X)).
\end{equation}
\end{example}
As we will need it later in this article let us finally prove the following
factorization property for $N_\order(\alpha_\mu)$.
\begin{lemma}
\label{NalphaFactorLem}
Let $\alpha \in \Gamma^\infty(T^*Q)$ and $t \in \mathbb R$. Then
\begin{equation} \label{NalphaFactor}
    N_\order (t\alpha)
    = \exp \left(-\im\order\lambda (\Delta + t\Fdiff{\alpha})\right)
    = \exp \left(t\Fdiff{\frac{\exp (- \im \order \lambda \Ds) - \id }
    {\Ds} \, \alpha}\right) \, \exp \left(-\im \order \lambda \Delta \right) .
\end{equation}
\end {lemma}
\begin {proof}
By a straightforward computation one finds
\[
    \frac{d}{dt} \exp\left( - \im \order \lambda (\Delta + t \Fdiff{\alpha})
    \right)
    = \exp \left(-\im \order \lambda (\Delta + t \Fdiff{\alpha}) \right)
    \left(\frac{\id - \exp \left( \ad \left( \im \order \lambda
    (\Delta + t\Fdiff{\alpha}) \right)\right)}
    {\ad \left(\Delta + t \Fdiff{\alpha}\right)}
    \left(\Fdiff{\alpha} \right)\right).
\]
Using (\ref {NPCom}) and the fact that $\Fdiff{\gamma}$
commutes with $\Fdiff{\gamma'}$ one concludes that the left-hand
side of (\ref {NalphaFactor}) satisfies the differential equation
\[
   \frac{d}{dt} \exp \left( - \im \order \lambda (\Delta + t \Fdiff{\alpha})
   \right) =
    \Fdiff{\frac{\exp \left(- \im\order \lambda \Ds \right) - \id}{\Ds}
    \, \alpha }
    \exp \left( - \im \order \lambda (\Delta + t\Fdiff{\alpha})\right).
\]
Now, the  right-hand side of (\ref{NalphaFactor}) solves this differential
equation with the correct initial condition. Hence the claim follows.
\end {proof}
\section {Star products for cotangent bundles with magnetic fields}
\label{MagFieldSec}
In this section we construct new star products from $\staro$ reflecting
the presence of a {\it magnetic field} on the
configuration space.
Hereby we will describe the magnetic field by a formal power series
$B \in \Gamma^\infty(\bigwedge^2 T^*Q)[[\lambda]]$ of \emph{closed}
two-forms on $Q$. If $B $ is exact, then
$A \in \Gamma^\infty(T^*Q)[[\lambda]]$ with $B =dA$ is called a
{\it vector potential} for $B$. Whenever $B$ is only closed but not exact
we shall speak of a {\it magnetic monopole}.
As before we only consider vector potentials and magnetic fields such that
(at least) in order $0$ the forms $A$ and $B$ are real.

Now let $\left\{ O_j \right\}_{j\in I}$ be an open cover of $Q$ by contractible
sets and $\left\{ A^{j} \right\}_{j \in I} $ a family of local vector
potentials for $B$, i.e.~$A^{j} \in \Gamma^\infty(T^*O_j)[[\lambda]]$
and $dA^{j} = B|_{O_j}$ for all $j \in I$.
Furthermore let $\Aj$ be the local time evolution operator which is
induced by $A^{j}$ on ${\mathcal C}^\infty (T^*O_j)[[\lambda]]$ at
time $t = 1$ according to Theorem \ref{EvolutionTheo}. Note that the operator
$\Aj$ depends on $\order$.
Now we define an associative product $\star_\order^j$ on
${\mathcal C}^\infty (T^*O_j)[[\lambda]]$ by
\begin{equation} \label {starswBDef}
    f \star_\order^j g :=
    \Aj \left( \big(\Aj\big)^{-1} (f) \staro \big(\Aj\big)^{-1} (g)
    \right), \qquad f,g \in {\mathcal C}^\infty (T^*O_j)[[\lambda]].
\end{equation}
In case $O_j \cap O_k$ is not empty, Theorem \ref{EvolutionTheo} (iii)
entails that $  \big({\mathcal A}^{k}\big)^{-1}  \Aj $
is the  evolution operator for $t =1$ of the closed one-form
$A^{j}|_{T^*(O_j \cap O_k)} - A^{k}|_{T^*(O_j \cap O_k)}  $. Hence by
Theorem \ref{EvolutionTheo} (iv) one concludes that
$\big({\mathcal A}^{k}\big)^{-1}  \Aj $ is an automorphism of
$ \left( {\mathcal C}^\infty (T^*(O_j \cap O_k))[[\lambda]] , \staro \right)$.
But then $ f \star_\order^j g = f \star_\order^k g$ follows for all
$f,g \in {\mathcal C}^\infty (T^*(O_j \cap O_k))[[\lambda]]$.
Therefore one can define a product $\starBo$ on
${\mathcal C}^\infty (T^*Q)[[\lambda]]$ by setting for all $j \in I$
\begin{equation}
   f \starBo g \big|_{T^*O_j} =
   f\big|_{T^*O_j} \star_\order^j g\big|_{T^*O_j}, \qquad
   f,g \in {\mathcal C}^\infty (T^*Q).
\end{equation}
The product $\starBo$ now does not depend on the
particular choice of the covering $\left\{O_j\right\}_{j \in I}$
and its local vector potentials $A^j$ but only on $B$,
as the above argument shows.
\begin{theorem}
\label{starBTheo}
Let $B \in \Gamma^\infty(\bigwedge^2 T^*Q)[[\lambda]]$ be a formal series of
closed two-forms $dB=0$. Then $\starBo$ has the following properties.
\begin {enumerate}
\item $\starBo$ is a star product on
      ${\mathcal C}^\infty (T^*Q)[[\lambda]]$ with respect to the symplectic
      form $\omega_{B_0} = \omega + \pi^*B_0$.
\item $\starsB$ is of standard order type.
\item One has $f \starBo g =
      N_\order^{-1} \left( N_\order (f) \starsB N_\order (g)\right)$ for all
      $f, g \in {\mathcal C}^\infty (T^*Q)[[\lambda]]$.
\item If $B = \lambda B_1$ then $\starBo$ is homogeneous. If
      $B = \cc B$ then $\starwB := \star_{1/2}^B$ satisfies
      $\overline{f \starwB g} = \overline{g} \starwB \overline{f}$
      for all $f,g$. If in addition $B$ has the shape $B=\sum_{r=0}^\infty
      (\im \lambda)^{2r} B_{2r}$ with real two-forms $B_{2r}$ then
      $\starwB$ is of Weyl type.
\end{enumerate}
\end{theorem}
\begin{proof}
It suffices to prove the theorem locally that means to show it under the
assumption $Q=O_j$. Now, in order $0$ in $\lambda$ the evolution
operator $\Aj$ is equal to the diffeomorphism
$\phi_j: \zeta_x \mapsto \zeta_x - A_0^{j}(x)$ of $T^*O_j$.
But as $\phi^*_j \omega = \omega + \pi^*B_0$, this entails (i).
Since part (ii) is a local statement it is sufficient to consider
functions having their support in some $\pi^{-1} (O_j)$. Then
$(\pi^* u) \star_0^B f = \Aj (((\Aj)^{-1} (\pi^*u)) \star_0 ((\Aj)^{-1} f))
= \Aj((\pi^*u) (\Aj)^{-1}f) = (\pi^*u) f$ since $\stars$ is of standard order
type and $\Aj \circ \pi^* = \pi^*$. Eq.~(\ref{DeveloperW}) entails (iii).
Due to theorem \ref {EvolutionTheo} (v) the homogeneity of
$\star_\order^{\lambda B_1}$ is obvious as well as the fact that complex
conjugation is an anti-automorphism of $\starwB$ for $B = \cc B$. In case
$B=\cc B$ and $B = \sum_{r=0}^\infty (\im \lambda)^{2r} B_{2r}$ the whole
construction of the star product $\starwB$ only involves the combination
$\im \lambda$ and all the differential and bidifferential operators
occuring as coefficients of the powers of this combination are real since
$\starw$ is of Weyl type implying that $\starwB$ is of Weyl type too.
\end{proof}

If $B $ is exact, then of course one can use any global one-form $A$ with
$B =dA$ to define $\Ao (1)$ and $\starBo$ directly   via Eq.~(\ref{starswBDef}).
The crucial point however lies in the fact that $\starBo$ can still be
constructed, if $B$ is only closed but not necessarily exact.

Let us now state an immediate consequence from the preceding theorem
concerning the star exponential which will ease some later calculations.
The star exponential $\Exp (tH)$ of a function
$H \in \mathcal C^\infty (T^*Q)[[\lambda]]$ with $t \in \mathbb R$
is defined as the unique solution of the differential equation
$\frac{d}{dt} \Exp (tH) = H * \Exp (tH)$ with initial condition
$\Exp (0) = 1$ (cf. \cite [Lem.~2.2]{BRW98a}).
\begin{corollary}
\label{CorStarExp}
\begin{enumerate}
\item
  For all $u,v \in {\mathcal C}^\infty (Q)[[\lambda]]$ one has
\begin{equation}
  (\pi^* u) \starBo (\pi^* v) = \pi^*(uv).
\end{equation}
\item
 The star exponential $\Exp (\pi^*u)$ with respect to $\starBo$ coincides
 with $\eu^{\pi^*u}$.
\end{enumerate}
\end{corollary}

Next we compute Deligne's characteristic  class of the star product $\starBo$.
This characteristic class has been introduced in \cite{Del95}
and classifies in a functorial way the isomorphy classes of star products on
a symplectic manifold $(M,\omega)$. It lies in the affine space
$-\frac{\im[\omega]}{\lambda} + H^2_{\mbox{\tiny dR}}(M)[[\lambda]]$ and can be
calculated by  methods of \v{C}ech cohomology.
Let us provide some details of the calculation as far as they are needed
for our purposes. For proofs and explicit arguments we refer
the reader to the exposition \cite{GR98}. At this instance we should
mention that our conventions differ from those used in \cite{GR98} by a
sign in the Poisson bracket and an additional factor $\im$ in the formal
parameter. If $\star$ is a star product on the symplectic manifold
$(M,\omega)$ there exists (cf.~\cite[Prop.~5.2]{GR98}) a good open cover
$\left\{ U_j  \right\}_{j\in J }$ of $M$ together with a family
$\left( D_j \right)_{ j \in J }$ of $\lambda$-{\it Euler derivations} of
$(\mathcal C^\infty (U_j),\star)$, i.e.~a family of derivations $D_j$ of
$\star$ over $U_j$ having the form
$$
   D_j = \lambda \frac{\partial}{\partial \lambda} + X_j  + D_j',
$$
where $X_j$ is conformally symplectic
$( \mathcal L_{X_j} \omega|_{U_j} =\omega_{U_j})$ and
$D_j' = \sum_{k\leq 1} \lambda^k D_{jk}'$ is a formal differential operator
over $U_j$.
As every $\lambda$-linear derivation is of the form
$\frac{1}{\lambda} {\rm ad}_\star (d)$ with $d \in \mathcal C^\infty (M)
[[\lambda]]$ there exist formal functions $d_{ij} \in \mathcal C^\infty
(U_i \cap U_j)[[\lambda]]$ fulfilling
$$
    D_j - D_i = \frac{1}{\lambda} {\rm ad}_\star (d_{ij})
$$
over $U_i \cap U_j$. Now, the sums $d_{ijk} = d_{ij} + d_{jk} + d_{ki}$
lie in $\C[[\lambda]]$ and define a $2$-cocycle whose \v{C}ech class
$d(\star) = [d_{ijk}] \in H^2 (M,\C)[[\lambda]]$ does not depend on the
choices made and is called {\it Deligne's intrinsic derivation-related class}.
\begin{definition} $({\it cf}.~\cite[{\it Def.~6.3}]{GR98})$
  The {\it characteristic class} $c(\star)$ of a star product $\star$ on
  $(M,\omega)$ is the element
  $c(\star) = -\frac{\im[\omega]}{\lambda} + \sum_{n=0}^\infty \, c(\star)^n
  \lambda^n $ of the affine space
  $-\frac{\im[\omega]}{\lambda} + H^2_{\mbox{\tiny dR}}(M)[[\lambda]]$ defined by
  \begin{equation}
  \nonumber
    c(\star)^0  = 2 (C_2^-)^{\#}, \quad
    \frac{\partial}{\partial \lambda} c(\star) (\lambda)  =
    \frac{1}{\lambda^2} d(\star).
  \end{equation}
  Hereby $(C_2^-)^{\#}$ is the two-form induced by the antisymmetric part
  $C_2^-$ of the bidifferential operator $C_2$ in the expansion of $\star$:
  $$
    (C_2^-)^{\#} (X_f ,X_g) =  C_2^- (f,g), \quad f,g \in \mathcal
    C^\infty(M)
  $$
  where $X_f$ denotes the Hamiltonian vector field with respect to
  $\omega$ that corresponds to $f$.
\end{definition}
In the particular case of the star products $\star_\order$ on $M = T^*Q$
the homogeneity operator $\mathsf H$ comprises a (global) $\lambda$-Euler
derivation for $\star_\order$, hence $d(\star_\order)=[0]$ follows.
Moreover, as $\starw$ is of Weyl type, the corresponding $ C_2$ is symmetric,
hence $(C_2^-)^{\#} =0 $ and
$c (\starw)= c(\star_\order) = [0]$ since $\omega$ is exact in this case.

So as a first step we compute the derivation related class to determine
$c(\star_\order^B)$. For that we need the following lemma. Note that it
involves the symmetric degree derivation $\degs$ which is defined by $\degs
\gamma = k\gamma$ for all $\gamma \in \Gamma^\infty (\bigvee^k T^*Q)$.
\begin{lemma}
\label{DeveloperHComLem}
Let $\gamma \in \Gamma^\infty (\bigvee T^*Q)[[\lambda]]$. Then
\begin{equation}
\label{HomFCom}
    [\mathsf H, \Fdiff{\gamma}] =
    \Fdiff {\left(\lambda\frac{\partial}{\partial\lambda} - \degs\right)
    \gamma}
\end{equation}
which for
$A = A_0 + \lambda A_1 + O(\lambda^2) \in \Gamma^\infty (T^*Q)[[\lambda]]$ entails
\begin{equation} \label{HomEvolutionCom}
    \Ao (t) \circ \mathsf H \circ \Ao(-t) =
    \mathsf H + t \Fdiff{
    \frac{\exp (\im\order\lambda\Ds) - \exp(-\im(1-\order)\lambda\Ds)}
         {\im\lambda\Ds} \left(\lambda\frac{\partial}{\partial \lambda} - \id
         \right)A}.
\end{equation}
\end{lemma}
\begin{proof}
Since (\ref{HomFCom}) is a relation between algebra morphisms and
derivations it is sufficient to prove it for the case where $\gamma$
is a function or one-form on $Q$, which is a simple computation.
But then (\ref{HomEvolutionCom}) follows
immediately from (\ref{DeveloperS}) and (\ref{HomFCom}) observing that $\phi_t^*
\circ \mathsf H \circ {\phi_{-t}}^* = \mathsf H - t \Fdiff{A_0}$
whith $\phi_t$ given as in (\ref{phitDef}).
\end{proof}

Now choose a good cover $\left\{ O_j  \right\}_{j\in J }$ of $Q$ and let the
$A^j$ and $\mathcal A^j$ like above. Since $\mathsf H$ is a global
$\lambda$-Euler derivation with respect to $\staro$ one obtains local
$\lambda$-Euler derivations $D_j$ with respect to the restriction of
$\starBo$ to $\mathcal C^\infty(T^*O_j)[[\lambda]]$ by defining $D_j:=
\mathcal A^j \circ \mathsf H \circ {\mathcal A^j}^{-1}$. From
(\ref{HomEvolutionCom}) it is obvious, that on $C^\infty(T^*(O_i
\cap O_j))[[\lambda]]$ one has
\[
D_j - D_i = \Fdiff{
    \frac{\exp (\im\order\lambda\Ds) - \exp(-\im(1-\order)\lambda\Ds)}
         {\im\lambda\Ds} \left(\lambda\frac{\partial}{\partial \lambda} - \id
         \right)(A^j - A^i)}.
\]
Now since $d(A^j - A^i) = 0$ on $O_j \cap O_i$ one can choose formal
smooth functions $c^{ij}$ on $O_j \cap O_i$ such that $d c^{ij} = A^j -
A^i$. From (\ref{adstaro}) we thus get $D_j - D_i = \frac{1}{\lambda}
\ad_{\staro} (-\im (\lambda\frac{\partial}{\partial\lambda} - \id )
\pi^* c^{ij})$ where we may, using the properties of $\mathsf F$, replace
$\ad_{\staro}$ by $\ad_{\starBo}$ since $\mathcal A^j \circ \pi^* = \pi^*$.
Now the combination $c^{ijk}= c^{ij} + c^{jk} + c^{ki}$ is
constant, hence induces a \v{C}ech cocycle whose cohomology class is equal
to the one of the magnetic field: $ [{c^{ijk}}]= [B]$ implying by the
above definition that $d (\starBo) = -\im(\lambda
\frac{\partial}{\partial \lambda}- \id )[\pi^*B]$, hence we have shown:
\begin{proposition}\label{DerivRelClassProp}
Let $B = B_0 + \lambda B_1 + O(\lambda^2)
\in \Gamma^\infty (\bigwedge^2 T^*Q)[[\lambda]]$ be a formal closed two-form
with real $B_0$. Then Deligne's intrinsic derivation-related class of $\starBo$
is given by $-\im(\lambda \frac{\partial}{\partial\lambda} - \id)
[\pi^* B]= \im [\pi^* B_0] -\im \sum_{n=2}^\infty (n-1)\lambda^n
[\pi^*B_n])$.
\end{proposition}
Again by definition of the characteristic class $c(\starBo)$ we get from the
preceeding proposition that $c(\starBo) = -\frac{\im}{\lambda}[\pi^*B -
\lambda \pi^*B_1 ] + c(\starBo)^0$, where $c(\starBo)^0$ has to be determined
as in the above definition. To this end we observe that by Theorem
\ref{starBTheo} the star products $\starBo$ for different $\order$ are
equivalent so that we may restrict our considerations to the
case $\order = 1/2$. A lengthy
but straightforward computation using the local equivalence transformation
$\mathcal A^j$ expanding it up to the first order $\lambda$ and the fact
that $\star_{1/2}$ is of Weyl type yields $C_2^-(f, g) = -\frac{\im}{2}
(\pi^*B_1)(X_f^{B_0}, X_g^{B_0})$, where the additional superscript $B_0$
indicates that we used the symplectic form $\omega_{B_0} = \omega + \pi^*B_0$ to
define the Hamiltonian vector fields. Obviously this result is independent
of the chosen local potential $A^j$ thus we get $c(\starBo)^0 =
-\im[\pi^* B_1]$. Consequently we obtain the following theorem.
\begin{theorem}
Let $B = B_0 + \lambda B_1 + O(\lambda^2)
\in \Gamma^\infty (\bigwedge^2 T^*Q)[[\lambda]]$ be a formal closed two-form
with real $B_0$. Then Deligne's characteristic class of $\starBo$ is given by
$-\frac{\im}{\lambda} [\pi^* B]$,
whence for all equivalence classes of star products for
$(T^*Q, \omega+ \pi^*B_0)$ there is a representative $\starBo$.
\end{theorem}
\begin{proof}
By the above considerations the assertion about Deligne's characteristic class
has already been shown. Since the de Rham cohomology of $T^*Q$ is
canonically isomorphic via $\pi^*$ to $H_{\mbox{\tiny \rm dR}} (Q)$
the final statement in the theorem holds as well.
\end{proof}

As a first application of our explicit formulas relating
$\starBo$ and $\staro$ we will prove finally in this section that all star
products $\starBo$ are strongly closed, that means integration over $T^*Q$
with respect to the corresponding volume form
$\Omega_{B_0} = \frac{(-1)^{[n/2]}}{n!} \omega_{B_0} \wedge \cdots \wedge
\omega_{B_0}$ is a trace functional for these star products
(see \cite {CFS92} for definitions). Note  that $\Omega_{B_0} = \Omega$
since $\pi^*{B_0}$ is a horizontal two-form.
\begin{lemma}
\label{AorderIntegralLem}
Let $A = A_0 + \lambda A_1 + O(\lambda^2) \in
\Gamma^\infty (T^*Q)[[\lambda]]$ be a formal one-form with
$A_0 = \cc A_0$ and let $\Ao$ be the corresponding time
development operator at $t = 1$. Then for all
$f \in \mathcal C^\infty_{\rm cpt} (T^*Q)[[\lambda]]$
\begin {equation} \label {AorderIntegral}
    \int_{T^*Q} f \Omega
    = \int_{T^*Q} \left( \Ao^{-1} f\right) \Omega.
\end{equation}
\end{lemma}
\begin{proof}
The fiber translating diffeomorphism induced by $A_0$ in $\Ao$ is
volume preserving since $\Omega_B = \Omega$. Thus we have only to
deal with the differential operator part of $\Ao$ as in
(\ref {DeveloperS}) which is the identity in lowest order of $\lambda$
and a sum of homogeneous differential operator of negative degree in
higher orders of $\lambda$. Due to \cite [Lem.~8.1] {BNW97b} these
higher orders do not contribute to the integral. This proves the
lemma.
\end{proof}
\begin{proposition}
Let $B = B_0 + \lambda B_1 + O(\lambda^2)
\in \Gamma^\infty (\bigwedge^2 T^*Q) [[\lambda]]$ be a closed
two-form with real $B_0$. Then for every $\order \in [0,1]$
the star product $\starBo$ is strongly closed.
\end{proposition}
\begin{proof}
We have to show that the integral over $T^*Q$ with respect to $\Omega$
of the commutator of two functions
$f, g \in \mathcal C^\infty_{\rm cpt} (T^*Q) [[\lambda]]$ vanishes. By a
partition of unity argument we may assume that the supports of
these functions are contained in a contractible open set $O$. Then
choose a local formal one-form $A$ such that $dA = B|O$. Due to
the definition of $\starBo$ and Lemma.~\ref {AorderIntegralLem} we
have
\[
    \int_{T^*Q} (f \starBo g - g \starBo f) \Omega
    =
    \int_{T^*O} \left(
    (\Ao^{-1} f) \staro (\Ao^{-1} g) -
    (\Ao^{-1} g) \staro (\Ao^{-1} f) \right) \Omega
    = 0
\]
since $\staro$ is homogeneous and thus strongly closed due to
\cite [Thm.~8.5] {BNW97b}. This proves the proposition.
\end {proof}
\section{Global representations for $\starBo$}
\label{RepMagFieSec}
After having constructed the product $\starBo$ of observables for the
case, where a magnetic field is present, we are now interested in
representations of these algebras.
It will turn out that, though the algebra structure depends only on $B$,
the representations have to be constructed by explicit use of the vector
potential $A$.
Interesting physical effects like the Aharonov-Bohm effect are only
visible by comparing the representations for varying vector potential $A$
but fixed $dA = B$.

Motivated by our investigations in \cite [Sect.~7] {BNW97b} we will
consider first the case where $A \in \Gamma^\infty(T^*Q)[[\lambda]]$ is
globally given with real $A_0$ and $B = dA$.
Then one can embed the graph
\begin{equation} \label {graphADef}
     L_{A_0} := \graph (A_0) =
     \left\{ \zeta_x \in T^*Q \; | \; \zeta_x = A_0 (x) \right\}
\end{equation}
into $T^*Q$ via $\iota_{A_0} : L_{A_0} \hookrightarrow T^*Q$.
Recall that $L_{A_0}$ is Lagrangian if and only if $d A_0 =0$.
Furthermore the relation
$\phi_{-t} (\iota(Q)) =\iota_{tA_0} (L_{tA_0})$ holds for all $t \in \R$.
Hence there exists an induced diffeomorphism
$\Phi_t: Q \to L_{tA_0}$ fulfilling
\begin{equation} \label {QPhiLADef}
    \phi_{-t} \circ \iota = \iota_{tA_0} \circ \Phi_t .
\end{equation}
We denote the pull-back $(\Phi_t^{-1})^* $ by
$U_t : {\mathcal C}^\infty (Q)[[\lambda]] \to {\mathcal C}^\infty
(L_{tA_0}) [[\lambda]]$ and clearly have $U_t^{-1} = \Phi_t^*$.
Using  $U = U_1$ we can now define a
representation ${\tilde\rep}^A_\order$ of the algebra
$\left( {\mathcal C}^\infty (T^*Q)[[\lambda]], \starBo\right)$ on the space
${\mathcal C}^\infty (L_{A_0})[[\lambda]]$ of formal wave functions on
$L_{A_0}$ by setting
\begin{equation} \label{tildeswrepADef}
    {\tilde\rep}^A_\order (f) u
    := \left(U \circ \repo \left({{\mathcal A}}^{-1} f\right)
    \circ U^{-1} \right) u, \qquad u \in  {\mathcal C}^\infty (L_{A_0}),
\end{equation}
where ${\mathcal A} = \Ao (1)$ is again the time development with
respect to $A$.
Obviously ${\tilde\rep}^A_\order$ is a representation with respect to
$\starBo$. Moreover, ${\tilde\rep}^A_\order$ is equivalent to the
representation $\repAo$ on $\mathcal C^\infty (Q)[[\lambda]]$ given by
\begin{equation} \label{swrepADef}
    \repAo (f) := \repo \left({{\mathcal A}}^{-1} f \right) .
\end{equation}
The operator $U$ then serves as an intertwiner between
${\tilde\rep}^A_\order$ and $\repAo$ and one has
$\repAo (f) = \srepA (N_\order f)$.

In order to compute $\srepA$ and thus the other representations
$\repAo$ and ${\tilde\rep}^A_\order$ explicitly we first prove
the following lemma.
\begin{lemma}
\label{DsLgammaLem}
  Let $\gamma \in \Gamma^\infty(\bigvee T^*Q)[[\lambda]]$ be a formal series of
  symmetric forms of degree $\ge 1$, and let $\gamma$ act on
  $\Gamma^\infty(\bigvee T^*Q)[[\lambda]]$ by left-multiplication.
  Then  for all $c \in \mathbb C[[\lambda]]$ and $t \in \mathbb R$
  the operator $\exp\left(c\Ds + t\gamma\right)$
  satisfies the factorization property
\begin{equation} \label {DsLgammaFactor}
  \exp\left(c\Ds + t \gamma\right) =
  \exp \left(t \frac{\exp(c\Ds) - \id}{c\Ds} \, \gamma\right)
  \, \exp (c\Ds).
\end{equation}
\end {lemma}
\begin {proof}
Since $c\Ds + t \gamma$ raises the symmetric degree at least by one,
the operator  $\exp\left(c\Ds + t\gamma\right)$ is well-defined as formal
series in the symmetric degree. By commutativity of $\vee$ and the
derivation property of $\Ds$ one has
$[c\Ds + t\gamma, \gamma'] = c\Ds\gamma'$. Now, proceeding analogously
as in the proof of Lemma \ref {NalphaFactorLem} one finds
\[
    \frac{d}{dt} \exp\left(c\Ds + t \gamma\right) =
    \left( \frac{\exp(c\Ds)-\id}{c\Ds} \, \gamma \right)
     \exp\left(c\Ds + t \gamma\right).
\]
The right-hand side of (\ref {DsLgammaFactor}) now solves this
differential equation, hence the claim follows.
\end {proof}

Next let us compute $\repAo$ explicitly for the case where
$A = \lambda A_1 + \lambda^2 A_2 + O(\lambda^3)$ starts in order one. This
condition simplifies computations since in this case no
diffeomorphism part is present. Nevertheless, it is the most important
situation for the following constructions.
\begin{theorem}
\label{repATheo}
Let
$A = \lambda A_1 + \lambda^2 A_2 + O(\lambda^3)
\in \Gamma^\infty(T^*Q)[[\lambda]]$
be a formal series of one-forms starting in order $\lambda$ and
$B = dA$ the corresponding magnetic field. Then one has for all
$f \in {\mathcal C}^\infty (T^*Q)[[\lambda]]$ and
$u \in {\mathcal C}^\infty (Q)[[\lambda]]$
\begin{equation} \label{wrepAFormel}
    \repAo (f) u = \sum_{l=0}^\infty
    \frac{(-\im\lambda)^l}{l!}
    \sum_{j_1, \ldots, j_l}
    \iota^*\left(\frac{\partial^l N_\order f}
    {\partial p_{j_1} \cdots \partial p_{j_l}}\right)
    \frac{1}{l!}
    \brac{\partial_{q^{j_1}} \otimes \cdots \otimes \partial_{q^{j_l}},
    \left( \Ds + \frac{\im}{\lambda} A\right)^l u } .
\end{equation}
\end {theorem}
\begin {proof}
Using the relation $\repAo (f) = \srepA (N_\order f)$ it suffices to show
Eq.~(\ref{wrepAFormel}) only for the case $\order=0$, so let us fix
$\order=0$ for the rest of the proof.
Now, the right-hand side of (\ref{wrepAFormel}) is independent of the
particular chosen canonical coordinate system, so it defines a global
object indeed. Moreover, both sides are sheaf morphisms in $u$, so it
is sufficient to show (\ref {wrepAFormel}) over the domain $U \subset Q$
of a chart on $Q$, and to assume
$f \in {\mathcal C}^\infty (T^*U)[[\lambda]]$ and
$u \in {\mathcal C}^\infty (U) [[\lambda]]$.
To achieve this consider first the space
\[
    {\mathcal T} (U) := \left( \prod_{k=0}^\infty {\mathcal C}^\infty (T^*U)
    \otimes_{{\mathcal C}^\infty (U)} \Gamma^\infty
    \left(\mbox{$\bigvee$}^k T^*U \right) \right)[[\lambda]].
\]
The space ${\mathcal T} (U)$ carries a natural commutative product $m$
given by the pointwise product of functions in the first and the
$\vee$-product in the second tensor factor.
Using canonical coordinates over $T^*U$ one can deform $m$ to
a product $\starloc$ by
\[
    F \starloc G := m \circ \exp \left( - \im \lambda \sum_j
    \partial_{p_j} \otimes  i_{\mbox{\tiny s}} (\partial_{q^j}) \right)
    (F \otimes G) ,
    \qquad F,G \in {\mathcal T} (U).
\]
Hereby $i_{\mbox{\tiny s}}  (\partial_{q^{i}}) G$ denotes the insertion of
$\partial_{q^{i}}$
in the symmetric part of $G$. Since $\partial_{p_j}$ and
$i_{\mbox{\tiny s}}(\partial_{q^j})$ are commuting
${\mathcal C}^\infty (U)$-linear derivations of
${\mathcal T} (U)$, $\starloc$ is indeed an associative product
deforming $m$. Let us denote by
$\sigma: {\mathcal T} (U) \to {\mathcal C}^\infty (T^*U)[[\lambda]]$ the
projection on the factor of symmetric degree zero.
Denoting the right-hand side of
(\ref{wrepAFormel}) shortly by $\pi (f) \, u$ we obtain due to Lemma
\ref{DsLgammaLem} and the particular form of $\starloc$
\[
\begin{split}
    \pi(f) \, u
    & = \iota^* \sigma \left(f \starloc
        \exp\left(\Ds + \frac{\im}{\lambda} A \right) u \right)
    =   \iota^*\sigma \left( f \starloc
        \left(\exp\left( \frac{\exp \Ds - \id}{\Ds}
        \frac{\im}{\lambda} A \right) \exp ( \Ds ) u \right) \right) \\
    & = \iota^*\sigma \left( f \starloc \exp \left(
        \frac{\exp \Ds - \id}{\Ds} \frac{\im}{\lambda} A \right)
        \starloc \exp ( \Ds)  u \right) \\
    & = \iota^*\sigma \left( \sigma \left( f \starloc \exp \left(
        \frac{\exp \Ds - \id}{\Ds} \frac{\im}{\lambda} A\right)
        \right) \starloc \exp ( \Ds ) u \right) .
\end{split}
\]
Now compute by Taylor expansion in the right factor
\[
    \sigma \left(f \starloc \frac{\exp \Ds - \id}{\Ds}
    \frac{\im}{\lambda} A \right) = \Fdiff{ \frac{\id -
    \exp \left(-\im \lambda\Ds \right)}{\im\lambda\Ds}
    A } f .
\]
Since
$\partial_{p_j} \frac{\exp \Ds-\id}{\Ds}\frac{\im}{\lambda} A = 0$ we
can use $\starloc$-products instead of $\vee$-products in the
exponential function. Hence, due to Theorem \ref{EvolutionTheo} and the
fact that $A$ starts in order $\lambda$ we obtain
\[
    \sigma \left(f \starloc \exp \left( \frac{\exp \Ds - \id}{\Ds}
    \frac{\im}{\lambda} A \right) \right) =
    \exp \left( \Fdiff{
    \frac{\id - \exp \left(-\im\lambda\Ds \right)}
    {\im\lambda\Ds} A} \right) f
    = {\mathcal A}_0^{-1} f.
\]
Using (\ref {SRepDef}) this finally leads to
\[
    \pi (f) \, u
    = \iota^*\sigma \left( {\mathcal A}_0^{-1} f \starloc \exp ( \Ds ) u
      \right)
    = \rep_0 ({\mathcal A}_0^{-1} f) u
    = \rep_0^A (f) u.
\]
\end {proof}

The interpretation of these formulas is that we replaced the covariant
derivative $\nabla_X$ acting as Lie derivative $\Lie_X$ on functions
$u \in {\mathcal C}^\infty (Q)$ by the covariant derivative
$\Lie_X + \frac{\im}{\lambda} A(X)$ where $A(X)$ acts as
left-multiplication on $u$. Thus we endowed the trivial
$\mathbb C$-bundle whose sections are the functions on $Q$ with a
non-trivial covariant derivative by adding a
\emph{formal connection one-form} $\frac{\im}{\lambda} A$ which
corresponds physically to the `minimal coupling' description and thus
justifies our interpretation of $A$ as `vector potential'.

Now we consider the Weyl ordered star product $\starwB$ for a real and
exact magnetic field $B = \cc B = dA$ with a real vector potential
$A \in \Gamma^\infty(T^*Q)[[\lambda]]$. In this case we already
mentioned that
$\starwB$ satisfies $\cc{f \starwB g} = \cc g \starwB \cc f$ for all
$f, g \in {\mathcal C}^\infty (T^*Q)[[\lambda]]$. Whence it is natural
to seek for positive $\mathbb C[[\lambda]]$-linear functionals on the
$^*$-algebra $\left({\mathcal C}^\infty (T^*Q)[[\lambda]], \starwB\right)$
such that the induced GNS representation is $\wrepA$ resp.~$\tilde \wrepA$.
Analogously to \cite[Sect.~7]{BNW97b} and guided by the
general idea of \cite[Prop.~5.1]{BNW97b} we denote the subspace of
those functions $f \in {\mathcal C}^\infty (T^*Q)$ such that
$\supp (\iota^*_{A_0} f)$ is compact in $L_{A_0}$ by
${\mathcal C}^\infty_{L_{A_0}} (T^*Q)$.
Since $\starwB$ is a differential star product,
${\mathcal C}^\infty_{L_{A_0}} (T^*Q)[[\lambda]]$ is a two-sided
ideal which is stable under complex conjugation. Denoting by $\Aw$ the
time evolution operator with respect to $A$ at time $t = 1$ the relation
\begin{equation}
\nonumber
    \Aw \left( {\mathcal C}^\infty_Q (T^*Q)[[\lambda]] \right)
    =
    {\mathcal C}^\infty_{L_{A_0}} (T^*Q) [[\lambda]]
\end{equation}
holds. Hence, the $\mathbb C[[\lambda]]$-linear functional
\begin{equation} % \label{omegaADef}
\nonumber
    \omega_A = \omega_\mu \circ \Aw^{-1} :
    {\mathcal C}^\infty_{L_{A_0}} (T^*Q)[[\lambda]] \to
    \mathbb C[[\lambda]]
\end{equation}
is well-defined. Moreover, it is positive with respect to $\starwB$,
since for all $f \in {\mathcal C}^\infty_{L_{A_0}} (T^*Q)[[\lambda]]$
\begin{equation}
\nonumber
    \omega_A (\cc f \starwB f)
    = \omega_\mu \left(\Aw^{-1} ( \cc f \starwB f)\right)
    = \omega_\mu \left(\cc{\Aw^{-1} f} \starw \Aw^{-1} f \right) \ge 0
\end{equation}
due to the reality of $\Aw$ and the definition of $\starwB$.
Explicitly, $\omega_A$ is given by
\begin{equation} \label{omegaAFormel}
    \omega_A (f) = \int_{L_{A_0}}
    \iota_{A_0}^* \exp \left( \Fdiff{
    \frac{\sinh \left(\frac{\im\lambda}{2}\Ds\right) }
    {\frac{\im\lambda}{2} \Ds} (A) - A_0} \right) (f) \; \mu_{A_0}
\end{equation}
according to (\ref{DeveloperW}) where the positive
density $\mu_{A_0} \in \Gamma^\infty(|\!\bigwedge^n\!|\, T^*L_{A_0})$
is given by $\mu_{A_0} = \left({\Phi^{-1}}\right)^*\mu$ with $\Phi$ as in
(\ref {QPhiLADef}) at $t = +1$. Clearly $\omega_A (\cc f \starwB f) = 0$
if and only if $\Aw^{-1} f \in \mathcal J_\mu$ is an element of the
Gel'fand ideal
$\mathcal J_\mu$ of $\omega_\mu$. Thus we have the following proposition
as slight generalization of \cite [Lem.~7.3] {BNW97b}. Note that due to
\cite [Cor.~1] {BW98a} the GNS representation extends to the whole
algebra.
\begin{proposition}
\label{GlobalGNSProp}
For real $A \in \Gamma^\infty(T^*Q)[[\lambda]]$ the Gel'fand ideal
$\mathcal J_A$ of $\omega_A$ is given by $\Aw(\mathcal J_\mu)$. The GNS
representation pre-Hilbert space
${\mathcal C}^\infty_{L_{A_0}} (T^*Q)[[\lambda]] \big/ \mathcal J_A$
is canonically isometric to
${\mathcal C}^\infty_{\rm cpt} (L_{A_0}) [[\lambda]]$
endowed with the $\C [[\lambda]]$-valued Hermitian product
\begin{equation}
    \langle u, v \rangle_A
    = \int_{L_{A_0}} \cc u \, v \, \mu_{A_0},
    \qquad
    u, v \in {\mathcal C}^\infty_{\rm cpt} (L_{A_0})[[\lambda]].
\end{equation}
A $\C [[\lambda]]$-linear unitary intertwiner is given by
$\psi_f \mapsto {\Phi^{-1}}^* \iota^* N_{1/2}(\alpha_\mu) \Aw^{-1} (f)$,
where $\psi_f$ denotes the equivalence class of
$f\in {\mathcal C}^\infty_{L_{A_0}} (T^*Q)[[\lambda]]$
in the GNS pre-Hilbert space.
The intertwiner has  inverse $v \mapsto \psi_{\pi^*\Phi^*v}$.
Moreover, the GNS representation induced by $\omega_A$ on
${\mathcal C}^\infty_{\rm cpt} (L_{A_0})[[\lambda]]$ coincides
with $\tilde \wrepA$ which  on
${\mathcal C}^\infty_{\rm cpt} (Q)[[\lambda]]$ is unitarily equivalent
to $\wrepA$.
\end{proposition}
\begin{corollary}
\label{PartIntDens}
  For all $\order \in [0,1]$ the scalar product $\brac{u, \repAo (f) v}_A$
  is given by
  \begin{equation}
    \brac{u, \repAo (f) v}_A
     = \int_{L_{A_0}} \cc u \, \repAo (f) v \, \mu_{A_0}
     = \int_{L_{A_0}} \cc
     {\repAo \left( N_{1-2\order} \cc f \right) u } \, v \, \mu_{A_0}
     = \brac{\repAo(N_{1-2\order} \cc f)u, v}_A.
  \end{equation}
\end{corollary}

As an  application of the above considerations and formulas we will briefly
examine the WKB expansion scheme as already discussed in
\cite[Sect.~7]{BNW97b} and \cite{BW97b}.
The starting point is a real-valued Hamiltonian function
$H \in {\mathcal C}^\infty (T^*Q)$ and a projectable Lagrangean submanifold
$\iota_{A_0}: L_{A_0} \hookrightarrow T^*Q$ determined by the graph of a
closed real one-form $A_0 \in \Gamma^\infty(T^*Q)$ such that the
Hamilton--Jacobi equation
\begin{equation} \label{HamiltonJacobi}
    H \circ A_0 = E
    \quad
    \mbox { resp. }
    \quad
    H|_{L_{A_0}} = E
\end{equation}
is satisfied for some energy value $E \in \mathbb R$. Then we have shown
in \cite [Thm.~7.4] {BNW97b} that the eigenvalue problem
\begin{equation} \label{WKBProblem}
    \wrep \left(\Aw^{-1} H\right) u = E u
\end{equation}
for an eigenfunction $u \in {\mathcal C}^\infty_{\rm cpt} (Q)[[\lambda]]$
or better an eigendistribution
$u \in {\mathcal C}^\infty_{\rm cpt} (Q)'[[\lambda]]$ results in the
usual transport equations of the WKB expansion by simple expansion of
(\ref {WKBProblem}) in powers of $\lambda$. Here $\Aw$ is the time
development operator corresponding to $A_0$ at time $t=+1$ which is in
this case an automorphism of $\starw$ since $dA_0 = 0$. In
\cite {BNW97b} we factored $\Aw^{-1}$ as $(\phi^{-1})^* \circ T_{-1}$
with some formal series of differential operators $T_{-1}$ which were
given by a recursion formula using iterated integrals. Now we can
determine $T_{-1}$ explicitly using Theorem \ref {EvolutionTheo}
resulting in
\begin{equation} \label{WKBEvoluter}
    \Aw^{-1} = \left( {\phi^{-1}} \right)^* \circ
       \exp \left(\Fdiff{ \frac{\sinh
       \left(\frac{\im\lambda}{2} \Ds\right)}
       {\frac{\im\lambda}{2}} (A_0) - A_0}\right) .
\end{equation}
Thus the expansion of $T_{-1}$ in powers of $\lambda$ is explicitly known
and can be used for a more explicit WKB expansion.
Note that the result of \cite [Lem.~7.5]{BNW97b} implying
$\Aw^{-1} H = (\phi^{-1})^* H$ for
$H \in {\mathcal P}(Q)[[\lambda]]$
at most quadratic in the momenta follows now directly form
(\ref {WKBEvoluter}) by `counting degrees'.
\section{Standard and Weyl representations on half-densities}
\label{RepHalDensSec}
As another application we will discuss how to carry over the
representations $\repo$ to representations on half-densities,
i.e.~on sections of the bundle
$\Gamma^\infty(|\!\bigwedge^n \!|^\half\, T^*Q)$.
Since the covariant derivative $\nabla$ acts on half-densities as well,
we define the standard order representation $\halfsrep$ of
${\mathcal C}^\infty (T^*Q)[[\lambda]]$ on
$\Gamma^\infty(|\!\bigwedge^n\!|^\half\, T^*Q)[[\lambda]]$ by
\begin{equation} \label {halfsrepDef}
    \halfsrep (f) \nu := \sum_{l=0}^\infty
    \frac{(-\im\lambda)^l}{l!}
    \sum_{j_1, \ldots, j_l}
    \iota^*\left(\frac{\partial^l f}
    {\partial p_{j_1} \cdots \partial p_{j_l}} \right) \frac{1}{l!}
    \brac{ \partial_{q^{j_1}} \otimes \cdots \otimes \partial_{q^{j_l}},
    \Ds^l \nu},
\end{equation}
where $\Ds^l \nu$ is viewed as a smooth section of
$\bigvee^l T^*Q \otimes |\!\bigwedge^n\!|^\half\, T^*Q$.
Clearly (\ref {halfsrepDef}) is globally defined, and we will
show in the sequel that this ad hoc quantization rule is in fact a
representation of a $\star_0^B$-product algebra for a particular
magnetic field $B$.

In order to compare $\halfsrep$ with $\srep$ we choose the square root
$\mu^\half$ of the positive density $\mu$ as a trivializing section for
the half-density bundle.
Note that $\nabla_X \mu^\half = \half \alpha_\mu (X)\, \mu^\half$
for all $X \in \Gamma^\infty (TQ)$.
Now one calculates easily
\begin{equation}
    \Ds (\gamma \otimes \mu^{1/2}) =
    \left(\Ds \gamma + \frac{1}{2} \, \alpha_\mu \vee \gamma\right) \otimes
    \mu^\half, \qquad \gamma \in \Gamma^\infty
    \big( \mbox{$\bigvee$}^l T^*Q\big).
\end{equation}
This implies the following result.
\begin{theorem}
Let $f \in {\mathcal C}^\infty (T^*Q)[[\lambda]]$ and
$\nu = v \mu^\half
\in\Gamma^\infty(|\!\bigwedge^n\!|^\half\, T^*Q)[[\lambda]]$
with $v \in {\mathcal C}^\infty (Q)[[\lambda]]$. If one  sets
$A = -\frac{\im\lambda}{2} \alpha_\mu \in \Gamma^\infty(T^*Q)[[\lambda]]$,
the relation
\begin{equation} \label{halfsrepA}
    \halfsrep (f) \nu = \left( \sum_{l=0}^\infty
    \frac{(-\im\lambda)^l}{l!}
    \sum_{j_1, \ldots, j_l}
    \iota^*
    \left(\frac{\partial^l f}
    {\partial p_{j_1} \cdots \partial p_{j_l}} \right)
    \frac{1}{l!}
    \brac{\partial_{q^{j_1}} \otimes \cdots \otimes \partial_{q^{j_l}},
    \left( \Ds + \frac{1}{2} \alpha_\mu \right)^l v}
    \right) \mu^\half
\end{equation}
holds, which implies
$\halfsrep (f) \nu = \left(\srepA (f) v\right) \mu^\half$.
\end{theorem}
Hence, the (non-canonical) trivialization of
$\Gamma^\infty(|\! \bigwedge^n \!|^\half\, T^*Q)[[\lambda]]$ induced by $\mu$
intertwines the (canonical) representation $\halfsrep$ of the
(canonical) star product $\starsB$ with the (non-canonical) representation
$\srepA$ associated to the purely imaginary vector potential $A$
and $B = dA$. Hereby {\it canonical} means depending only on $\nabla$,
whereas {\it non-canonical} means depending on $\mu$ as well. Note that by
definition of $\alpha_\mu$ the relation $\tr R = - d\alpha_\mu$ is true,
where $\tr R$ denotes the trace of the curvature tensor $R$ of $\nabla$
(cf.~\cite [Eqn.~(21)]{BNW97b}). So the `magnetic field' fulfills
\begin{equation}
      B = \frac{\im\lambda}{2} \tr R,
\end{equation}
hence $\star_0^B$ does not depend on the particular choice of $\mu$ indeed,
but only on $\nabla$. Thus we will denote $\star_0^B$ for this particular $B$
also by $\halfstars$.
\begin{corollary}
For $B = \frac{\im\lambda}{2} \tr R$ the star product $\halfstars = \starsB$
has a canonical representation of standard order type on half-densities
$\Gamma^\infty(|\! \bigwedge^n \!|^\half\, T^*Q)[[\lambda]]$. Explicitly it is
given by (\ref {halfsrepA}). Moreover, $\halfstars$ is homogeneous and
equivalent to $\stars$ via the time evolution operator $\mathcal A_0 (1)$
corresponding to $A = -\frac{\im\lambda}{2} \alpha_\mu$. The two star
products $\stars$ and $\halfstars$ of standard order type coincide if and only
if the connection $\nabla$ is unimodular, i.e.~if $\tr R = 0$. In this case
$\halfsrep$ is a representation of $\stars$ on half-densities.
\end {corollary}
Note that in general $\stars \ne \halfstars$ unless $\nabla$ is unimodular.
Nevertheless the unimodular case is the most important one, since e.g.~any
Levi-Civita connection to a given metric on $Q$ is unimodular.

The space of complex-valued smooth half-densities with compact support
$\Gamma^\infty_{\rm cpt} (|\!\bigwedge^n\!|^\half\, T^*Q)$
carries a natural pre-Hilbert space structure via the Hermitian product
\begin{equation} \label{HalfProd}
    \brac{\nu,\rho} = \int_Q \cc \nu  \rho \, , \qquad
    \nu,\rho \in
    \Gamma^\infty_{\rm cpt} (\mbox{$|\!\bigwedge^n\!|$}^\half\, T^*Q).
\end{equation}
It naturally extends to a $\mathbb C[[\lambda]]$-valued
Hermitian product on
$\Gamma^\infty_{\rm cpt} (|\!\bigwedge^n\!|^\half\, T^*Q)[[\lambda]]$.
Therefore it is reasonable to ask for the formal adjoint of $\halfsrep (f)$.
One can now calculate $\langle \halfsrep (f) \nu, \rho \rangle$
either by  intrinsic partial integration analogously to
\cite[Lem.~4.1]{BNW97b} or by using a particular trivialization as we
shall do in the following.
\begin{proposition}
For $f \in {\mathcal C}^\infty (T^*Q)[[\lambda]]$ one has
\begin{equation} \label{NalphaAsCC}
    \left( N_{1/2} (\alpha_\mu) \right)^2 \cc{{\mathcal A}_0^{-1} f} =
    {\mathcal A}_0^{-1} \left( N_{1/2} (0) \right)^2 \cc f,
\end{equation}
and for any two half-densities
$\nu, \rho \in
\Gamma^\infty_{\rm cpt}(|\!\bigwedge^n\!|^\half\, T^*Q)[[\lambda]]$
one has
\begin{equation} \label{halfsrepAdjoint}
    \left\langle \nu, \halfsrep (f) \rho \right\rangle =
    \left\langle \halfsrep \left( \left( N_{1/2}(0)\right)^2 \cc f \right)
    \nu, \rho \right\rangle .
\end{equation}
\end {proposition}
\begin {proof}
The first part is a straightforward computation using Lemma
\ref {NalphaFactorLem} and (\ref {NPCom}). Then the second part follows
from \cite [Lem.~4.1]{BNW97b} and (\ref {NalphaAsCC}) using the
trivialization induced by $\mu$.
\end {proof}

Thus $\halfsrep$ is not a unitary representation, but due to
(\ref {halfsrepAdjoint}) we can define a coresponding Weyl star product
$\halfstarw$ together with a representation  by setting
\begin{equation} \label{halfWeylDef}
    f \halfstarw g = \big( N_{1/2} (0) \big)^{-1}
    \left( \big( N_{1/2} (0)  f \big) \halfstars
    \big( N_{1/2} (0 ) g \big) \right)
\end{equation}
and $\halfwrep (f) = \halfsrep \big( N_{1/2} (0) f \big)$.
Hence we obtain the formula
\begin{equation} \label{halfwrepFormel}
    \halfwrep (f) \nu = \sum_{l=0}^\infty
    \frac{(-\im\lambda)^l}{l!}
    \sum_{j_1, \ldots, j_l}
    \iota^*\left(\frac{\partial^l N_{1/2} (0) f}
    {\partial p_{j_1} \cdots \partial p_{j_l}} \right) \frac{1}{l!}
    \brac{\partial_{q^{j_1}} \otimes \cdots \otimes \partial_{q^{j_l}},
    \Ds^l \nu} .
\end{equation}
Thus $\halfwrep$ comprises a representation with respect to $\halfstarw$
on the space of half-densities such that
\begin{equation} \label{halfwrepAdjoint}
    \left\langle \nu, \halfwrep (f) \rho \right\rangle=
    \left\langle \halfwrep (\cc f)\nu, \rho \right\rangle,
    \qquad \nu, \rho \in
    \Gamma^\infty_{\rm cpt} (\mbox{$|\!\bigwedge^n\!|$}^\half\,
    T^*Q)[[\lambda]].
\end{equation}
Note that $\halfstarw$ is defined by use of the connection $\nabla$ alone and
does not depend on the choice of $\mu$.
\begin{proposition} \label{halfWeylProp}
  The star product $\halfstarw$ is a homogeneous star product of Weyl
  type, and
  \begin{equation} \label {CChalfstarw}
    \cc{f \halfstarw g} = \cc g \halfstarw \cc f,
    \qquad f, g \in {\mathcal C}^\infty (T^*Q)[[\lambda]].
  \end{equation}
  Moreover, $\halfstarw$ is equivalent to the product $\starwB$ with
  $B = \frac{\im\lambda}{2}\tr R$ as well as to the product  $\starw$.
  An algebra isomorphism from $\starwB$ to $\halfstarw$ is given by
  $ \left( N_{1/2} (\alpha_\mu) \right)^{-1} N_{1/2} (0)$, one from
  $\starw$ to  $\halfstarw$ by
\begin{equation} \label{BDef}
    \mathcal B_\mu
    := \Aw^{-1}\left( N_{1/2}(\alpha_\mu)\right)^{-1} N_{1/2}(0)
    = \exp \left(\Fdiff{
    \frac{\cosh \left(\frac{\im\lambda}{2} \Ds\right) - \id}
    {\Ds} \alpha_\mu}\right) .
\end{equation}
\end {proposition}
\begin {proof}
The equivalence of $\halfstarw$ and $\starwB$ resp.~$\starw$ can be
shown by a straightforward computation using the definitions.
The explicit formula for $\mathcal B_\mu$ follows from
Lemma \ref{NalphaFactorLem}, the fact that
$\cc N_{1/2} (\alpha_\mu) = \left( N_{1/2}(\alpha_\mu)\right)^{-1}$, and
the explicit form of $\Aw$ according to (\ref{DeveloperS}).
Finally (\ref{CChalfstarw}) follows from (\ref{BDef}) and (\ref{CCAnti})
implying that $\halfstarw$ is of Weyl type since the whole construction
depends only on combination $\im\lambda$.
\end {proof}
\begin{remark}
If the connection $\nabla$ is unimodular, e.g.~the Levi-Civita connection
of a Riemannian metric on $Q$, we always have $\starwB = \starw$ for any
choice of $\mu$. Now one can ask for a covariantly constant density
$\mu_0$ or in other words for a density $\mu_0$ fulfilling
$\alpha_{\mu_0} = 0$. Note that for a non unimodular $\nabla$ such a
density cannot exist. Now, if $\nabla\mu_0 = 0$ and $Q$ is connected,
then $\mu_0$ is uniquely determined up to a (positive) scalar multiple.
Making the Ansatz $\mu_0 = \eu^{-\varphi}\mu$ in the unimodular case we
find that $\mu_0$ is covariantly constant if and only if $\alpha_\mu$
is not only closed but exact with $\alpha_\mu = d\varphi$. In particular we
succeed in finding a covariantly constant density for any given unimodular
connection in case $H^1_{\mbox{\rm\tiny dR}} (Q) = \{0\}$.
Now having found such a density one has even $\halfstarw = \starw$ and
thus $\halfwrep$ is a representation for $\starw$. Note furthermore that
even in this case $\halfstarw$ does  in general \emph{not} coincide with
the homogeneous Fedosov star product $\starf$ of Weyl type constructed in
\cite[Sect.~3]{BNW97a}.
\end{remark}

Finally we would like to show that the representation $\halfwrep$ is a
GNS representation.
Let $\nu \in \Gamma^\infty(|\! \bigwedge^n\!|^\half\, T^*Q)$ be
an arbitrary nowhere vanishing half-density and $\mu = \cc \nu \nu$
its square which comprises a positive density.
Define $\alpha_\mu$, etc.~with respect to this $\mu$.
Then consider the  $\mathbb C[[\lambda]]$-linear functional
\begin{equation} \label{omeganuDef}
    \omega_\nu : \: {\mathcal C}^\infty_Q (T^*Q)[[\lambda]] \rightarrow
    \mathbb C[[\lambda]],  \:
    f \mapsto \omega_\nu (f) = \int_Q \iota^* \mathcal B_\mu (f) \, \mu.
\end{equation}
Clearly it is  well-defined and local since $\mathcal B_\mu$ is a formal
series of differential operators.
Moreover, $\omega_\nu (f) = \omega_\mu (\mathcal B_\mu (f))$,
which allows us to apply the general result on the `pull back' of
GNS representations \cite [Prop.~5.1]{BNW97b} since due to
Proposition \ref{halfWeylProp} the map $\mathcal B_\mu$
is a \emph{real} algebra isomorphism between $\halfstarw$ and $\starw$.
Note that, though $\omega_\nu$ is only defined on a twosided ideal,
the GNS representation again extends to the whole algebra due to
\cite[Cor.~1] {BW98a}.
\begin {proposition}
Let $\nu \in \Gamma^\infty(|\!\bigwedge^n\!|^\half\, T^*Q)$ be a nowhere
vanishing half-density.
Then $\omega_\nu$ is a positive $\mathbb C[[\lambda]]$-linear functional
with Gel'fand ideal
\begin{equation} \label {omeganuGIdeal}
  \mathcal J_\nu = \left\{f \in {\mathcal C}^\infty_Q (T^*Q)[[\lambda]]
  \; \big| \; \iota^* N_{1/2}(\alpha_\mu) {\mathcal B}_\mu (f) = 0 \right\}.
\end{equation}
The GNS pre-Hilbert space
${\mathcal C}^\infty_Q (T^*Q)[[\lambda]] \big/ \mathcal J_\nu$ is canonically
isometric to
$\Gamma^\infty_{\rm cpt} (|\!\bigwedge^n\!|^\half\, T^*Q)[[\lambda]]$ via
\begin{equation} \label{nuIso}
    \psi_f \mapsto
    \left(\iota^* N_{1/2}(\alpha_\mu) {\mathcal B_\mu} f \right) \nu
    \quad
    \mbox { and its inverse }
    \quad
    v \nu \mapsto \psi_{\pi^* v},
\end{equation}
where $f \in {\mathcal C}^\infty_Q (T^*Q)[[\lambda]]$ and
$v \in {\mathcal C}^\infty_{\rm cpt} (Q)[[\lambda]]$.
The induced GNS representation of ${\mathcal C}^\infty (T^*Q)[[\lambda]]$ on
$\Gamma_{\rm cpt}^\infty (|\!\bigwedge^n\!|^\half\, T^*Q)[[\lambda]]$
coincides with $\halfwrep$.
\end {proposition}
\begin {proof}
This can be either checked explicitly by using the defining formula for
$\mathcal B_\mu$ or by using \cite[Prop.~5.1]{BNW97b}.
\end {proof}
\section{Aharonov-Bohm representations}
\label{ABSec}
Now we shall consider the question of equivalence classes of representations
$\repAo$ for different vector potentials but fixed magnetic field
$B = dA $. So let $A, A' \in \Gamma^\infty(T^*Q)[[\lambda]]$ both fulfill
$B = dA = dA' $ but not necessarily $A = A'$. Then we call the representations
$\repAo$ and $\rep_\order^{A'}$ \emph{equivalent}, if there exists a bijective
$\mathbb C[[\lambda]]$-linear \emph{intertwining} operator
$U: {\mathcal C}^\infty (Q)[[\lambda]] \to {\mathcal C}^\infty
(Q)[[\lambda]]$, i.e.~$U$ satisfies for all
$f \in {\mathcal C}^\infty (T^*Q)[[\lambda]]$ the relation
\begin{equation} \label {RepEquiDef}
    U^{-1} \rep_\order^{A'} (f) \, U = \repAo (f).
\end{equation}
In the case, where $A$ and $ A'$ are real and $\starBo$
is the Weyl product $\starwB$,
the representations $ \wrep^{A}$ and $\wrep^{A'}$
are called \emph{unitarily equivalent}, if in addition $U$ is unitary with
respect to the Hermitian product of Prop.~\ref {GlobalGNSProp}.
\begin {lemma}
\label {ULeftMultLem}
\label {NiceEquiLem}
An operator $U: {\mathcal C}^\infty (Q)[[\lambda]] \to {\mathcal C}^\infty
(Q)[[\lambda]]$ is an intertwiner for $\repAo$ and $\rep_\order^{A'}$
if and only if for all $f \in {\mathcal C}^\infty (Q)$
\begin{equation}
  U^{-1} \srep (f) \, U = \srep^{A-A'} (f).
\end{equation}
In that case $U$ is the left multiplication by a function in
${\mathcal C}^\infty (Q)[[\lambda]]$ also denoted by $U$. Moreover $U$
is unique up to an invertible factor in $\mathbb C[[\lambda]]$.
\end{lemma}
\begin {proof}
The first part of the lemma is a simple consequence of (\ref {AswAsw}),
(\ref{tildeswrepADef}) and (\ref{swrepADef}).
For the second part check that $U$  commutes with all left multiplications,
i.e.~$U u v = u U v$ for all $u,v \in {\mathcal C}^\infty (Q) [[\lambda]]$.
Hence $U$ is a left multiplication by a function in
${\mathcal C}^\infty (Q) [[\lambda]] $ as well.
\end {proof}

Due to this lemma the equivalence of two representations
$\repAo$ and  $\rep_\order^{A'}$ for a fixed value of the ordering
parameter $\order$ entails equivalence of $\rep^A_{\order'}$ and
$\rep^{A'}_{\order'}$ for any other value $\order'$ as well.
Hence we will restrict our considerations for a moment only to the
case where $\order = \frac 12$ and $B = 0$. The following lemma then
is immediately checked using the fact that the corresponding time
development operator is an \emph{inner} automorphism, see
e.g.~\cite [Lem.~2.3]{BRW98a} and Corollary \ref{CorStarExp}.
\begin{lemma}
\label{ExactEquiLem}
Let
$A = \lambda A_1 + \lambda^2 A_2 + O(\lambda^3)
\in \Gamma^\infty(T^*Q)[[\lambda]]$
be an exact one-form with $A = dS$ and
$S = \lambda S_1 + \lambda^2 S_2 + O(\lambda^3)
\in {\mathcal C}^\infty (Q)[[\lambda]]$. Then
$\wrep$ and $\wrepA$ are equivalent, where an intertwiner
$U: \mathcal C^\infty (Q)[[\lambda]] \to \mathcal C^\infty (Q)[[\lambda]]$
is given by the left multiplication with
$\exp \left(\frac{\im}{\lambda} S\right)$. If $A$ is real and
$S$ has also been chosen to be real, then the restriction of $U$
to ${\mathcal C}^\infty_{\rm cpt} (Q)[[\lambda]]$ is unitary.
\end {lemma}

Now choose a good cover $\left\{ O_j \right\}_{j \in I}$ of $Q$,
i.e.~assume that all finite intersections of the $O_j$ are contractible.
Next let $\left\{ S^j \right\}_{j \in I}$
be a formal \v Cech cochain integrating
$A$, or in other words a family of formal functions
$S^j = \lambda S^j_1 + \lambda^2 S^j_2 + O(\lambda^3)
\in {\mathcal C}^\infty (O_j)[[\lambda]]$ fulfilling
$d S^j  = A|_{O_j}$ for all $j \in I$.
Next let $U$ be an intertwiner from $\wrep$ to $\wrepA$,
and choose local logarithms $T^j \in \mathcal C^\infty (O_j) [[\lambda]]$
for $U$ that means $U\big|_{O_j} = \exp \left( 2 \pi \im \, T^j\right)$.
Then $(T^j - T^k)\big|_{O_j \cap O_k} \in \Z$ holds for all nonempty
intersections $O_j \cap O_k$.
As the center of $\starw$ consists only of the constant functions,
one can find by the preceding lemma $c^j \in \C$ such that
$2 \pi T^j = \frac{S^j}{\lambda} + c^j $. Therefore
we may assume without restriction that
\begin{equation}
\label{Integralitaet}
\frac{(S^j - S^k)}{\lambda} \in 2 \pi \Z
\end{equation}
for all pairs $j,k$ with $O_j \cap O_k \neq \emptyset$.
Let us at this point abbreviate the notation by calling a formal closed
$k$-form $\alpha = \lambda \alpha_1 + \lambda^2 \alpha_2 + O(\lambda^3)
\in \Gamma^\infty (\bigwedge^k T^*Q)[[\lambda]]$ {\it integral}, if it
satisfies the following condition
\begin{enumerate}[$(${\rm IC}$)$]
\item
    \mbox{ }\hspace{2mm} $\frac{1}{2\pi} \alpha_1$ defines an integral
    cohomology class and $\alpha_l$ is exact for all $l \geq 2$.
\end{enumerate}
By the above the existence of an intertwiner $U$ then implies that
(IC) holds for $A$.

Vice versa suppose now that $A$ is integral. Then one can find
$S^j$ such that (\ref{Integralitaet}) holds again. Consequently one can define
$U \in {\mathcal C}^\infty (Q)[[\lambda]]$ by requiring
\begin{equation}
\label{UglobalDef}
   U\big|_{O_j} = \exp \left( \frac{\im}{\lambda} S^j \right) =
   \iota^* \Exp \left( \frac{\im}{\lambda} \pi^*S^j \right)
\end{equation}
for all $O_j$. By Corollary \ref{CorStarExp} and Lemma \ref{ExactEquiLem}
the operator $U$ then is an intertwiner from $\wrep$ to $\wrepA$.
Altogether we thus obtain the following theorem which completely describes
the equivalence classes of representations.
\begin {theorem}
Let $A, A' \in \Gamma^\infty(T^*Q)[[\lambda]]$ be vector potentials of the
form $A = \lambda A_1 + \lambda^2 A_2 + O(\lambda^3)$ resp.
$A' = \lambda A'_1 + \lambda^2 A'_2 + O(\lambda^3)$ fulfilling $B = dA =  dA'$.
Then the $\starBo$-representations $\repAo$ and $\rep_\order^{A'}$ are
equivalent if and only if the difference $A-A'$ is integral,
whence the equivalence classes are parametrized by
$\lambda H^1_{\mbox{\tiny \rm dR}} (Q ) / H^1_{\mbox{\tiny \rm dR}} (Q,\Z ) +
\lambda^2 H^1_{\mbox{\tiny \rm dR}} (Q )[[\lambda]] $.
In case, where in addition $A$ and $A'$ are real,
the representations $\wrepA$ and $\wrep^{A'}$ are even unitarily equivalent
with an intertwiner $U$ given locally by
$U|_{O_j} = \exp (\frac{\im}{\lambda} S^j)$, where $\{S^j\}_{j \in I}$ is
a real \v{C}ech cochain integrating $A-A'$.
\end {theorem}
In the case  $B=0$ we have a canonical reference
representation namely $\repo$ which corresponds to $A = 0$. Hence we
have also a canonical equivalence class of representations. Considering
the particular case $A = \lambda A_1$ without higher order terms the resulting
equivalence classes are parametrized by $H^1_{\mbox{\tiny dR}} (Q)$ modulo
integral classes. We shall call a representation $\rep_\order^{\lambda A_1}$ an
\emph{Aharonov-Bohm} representation if it is not equivalent to the
{\it vacuum} representation $\repo$ since in this case Aharonov-Bohm like
effects are possible. Thus the well-known integrality condition for the
vector potential can be found in deformation quantization as well.
But note that the corresponding star product $\starBo$ and thus the
algebra of observables is unaffected.
\section{Magnetic monopoles and non-trivial bundles}
\label{MagMonoSec}
In this section we will investigate the representation theory for
quantized algebras with product $\starBo$, where the magnetic field $B$ is
no longer exact but only closed. More precisely we will show under which
conditions representations of $\starBo$ can be constructed which
correspond to the  Schr\"{o}dinger-like representations $\repo$.
Throughout the whole section we consider only the case where
$B = \lambda B_1 + \lambda^2 B_2 + O(\lambda^3)$ starts in order
$\lambda$.

First of all we should mention that it is possible to construct
\emph{local} representations for all $\starBo$-algebras in the following
way: let $O \subset Q$ be a contractible open subset of $Q$ and choose
by Poincar\'{e}'s Lemma a local vector potential
$A \in \Gamma^\infty (T^*O)[[\lambda]]$ fulfilling $B|_O=dA$.
Then we can set
\begin{equation}\label {LocalRepDef}
    \rep_\order^O (f) u = \repAo (f|_{T^*O}) u, \qquad
  f \in \mathcal C^\infty (T^*Q)[[\lambda]], \:
  u \in \mathcal C^\infty (O)[[\lambda]],
\end{equation}
and thus obtain a $\starBo$-representation of the whole algebra on
$\mathcal C^\infty (O)[[\lambda]]$.
Definitely, these representations are not quite satisfactory
since they do not reflect the global nature of $B$, but by gluing them
together appropriately we will obtain the global representations we are
looking for.
The following lemma which is proved by straightforward computation
using Corollary \ref{CorStarExp} will
provide the essential tool for the construction.
\begin{lemma}
  Let $B = \lambda B_1 + \lambda^2 B_2 + O(\lambda^3)
  \in \Gamma^\infty (\bigwedge^2 T^*Q)[[\lambda]]$ be a closed two-form and
  $O_1, O_2 \subset Q$ two open subsets such that
  $O_1, O_2$ and the intersection $O_1 \cap O_2$ are contractible.
  Moreover, let
  $A^{(i)} = \lambda A^{(i)}_1 + \lambda^2 A^{(i)}_2 + O(\lambda^3) \in
  \Gamma^\infty (T^*O_i)[[\lambda]]$, $i = 1,2$ be local vector potentials
  for $B$. Finally choose
  $S^{(12)} = \lambda S^{(12)}_1 + \lambda^2 S^{(12)}_2 + O(\lambda^3) \in
  \mathcal C^\infty (O_1 \cap O_2)[[\lambda]]$ such that
  $dS^{(12)} = (A^{(1)} - A^{(2)})|_{O_1 \cap O_2}$.
  Then  the relation
  \begin{equation} \label{LocalRepEqui}
    \rep_  \order^{O_1} (f) u =
    \eu^{-\frac{\im}{\lambda} S^{(12)}} \rep_\order^{O_2} (f)
    \, \eu^{\frac{\im}{\lambda} S^{(12)}} u
  \end{equation}
  is fulfilled for all
  $f \in \mathcal C^\infty (T^*Q)[[\lambda]]$ and
  $u \in \mathcal C^\infty(O_1 \cap O_2)[[\lambda]]$.
\end{lemma}
The lemma now suggests that one should view the local functions $u$ as
coefficient functions for a particular {\it local trivialization} of a
{\it $\mathbb C[[\lambda]]$-line bundle} having the
$\eu^{-\frac{\im}{\lambda} S^{(ij)}}$ as transition functions.
We will not further develop such a theory of formal line bundles though this
could be done easily using standard sheaf cohomological methods extended to
modules over $\mathbb C[[\lambda]]$. Rather, we will construct an ordinary
complex line bundle out of the given data together with a formal
connection under the assumption that the magnetic field satisfies the
integrality condition (IC), see e.g.~\cite [p.~97--107] {Wells80}.
Then there exists a complex line bundle $\pi_L: L \rightarrow Q$
unique up to equivalence with Chern class $[\frac{1}{2\pi} B_1]$.
Moreover, we can find a connection $\nabla^L$ on $L$ such
that the curvature of $\nabla^L$ coincides with $\im B_1$ and
vector potentials for all $B_k$ with $k\geq 2$, i.e.~$B_k = dA_k$.
Let $A = \sum_{k=2}^\infty \lambda^k A_k$ then the mapping
\begin{equation}
  \nabla^A : \Gamma^\infty(L) [[\lambda]]
  \rightarrow
  \Gamma^\infty (T^*Q \otimes L)[[\lambda]] = \Omega^1(L)[[\lambda]]
  \quad s \mapsto
  \left(\nabla^L + \frac{\im}{\lambda} A \right) s
\end{equation}
comprises a formal connection on $L$ which we call {\it adapted} to $B$.
The equivalence classes of such $\nabla^A$ are parametrized by
$\lambda H_{\mbox{\tiny dR}}^1(Q) / H^1_{\mbox{\tiny dR}}(Q,\Z)$ $+
\lambda^2 H_{\mbox{\tiny dR}}^1(Q)[[\lambda]]$.
As usual let $\DL$ be the corresponding formal symmetrized covariant
derivative defined as in (\ref {SymCovDef}) using $\nabla^A$.
\begin {theorem}
Let $B = \lambda B_1 + \lambda^2 B_2 + O(\lambda^3) \in
\Gamma^\infty (\bigwedge^2 T^*Q)[[\lambda]]$ be a closed two-form
which satisfies the integrality condition $(${\rm IC}$)$.
Furthermore, let $\pi_L : L \to Q$ be a
complex line bundle with Chern class $[\frac{1}{2 \pi} B_1]$
and $A = \lambda^2 A_2 + O(\lambda^3) \in \Gamma^\infty (T^*Q)[[\lambda]]$
such that $dA = B - \lambda B_1$ whence $\nabla^A$ is a formal connection
adapted to $B$. Denote for every
$f \in {\mathcal C}^\infty (T^*Q) [[\lambda]]$ by
$\repLAo (f) : \Gamma^\infty (L) [[\lambda]] \to
\Gamma^\infty (L) [[\lambda]]$ the operator fulfilling
\begin{equation} \label{LwrepDef}
    \repLAo (f) s = \sum_{l=0}^\infty
    \frac{(-\im\lambda)^l}{l!}
    \sum_{j_1, \ldots, j_l}
    \iota^*\left(
    \frac{\partial^l N_\order f}{\partial p_{j_1} \cdots \partial p_{j_l}}
    \right) \frac{1}{l!}
    \brac{ \partial_{q^{j_1}} \otimes \cdots \otimes \partial_{q^{j_l}},
    \left(\DL\right)^l s}, \qquad
    s \in \Gamma^\infty (L)[[\lambda]].
\end{equation}
Then $\repLAo$ comprises a representation of $\starBo$ on
$\Gamma^\infty (L)[[\lambda]]$. Using a local trivializing section
$e_O : O \to L$, one has the local representation
\begin{equation} \label{LRepLocal}
    \left.\repLAo(f) s \right|_{O} =
    \left(\rep_\order^O  (f|_{T^*O}) u \right) e_O,
\end{equation}
where $s|_O = u \, e_O$ with $u \in \Gamma^\infty (O)[[\lambda]]$ unique
and where the local formal one-form $A'$ used for $\rep_\order^O$ as
in (\ref {LocalRepDef}) is determined by
$\nabla^A_X e_O = \frac{\im}{\lambda} A'(X) e_O$, i.e.~the local vector
potential satisfies $A' = \lambda A_1 + A$ with $d A_1 = B_1|_O$.
\end {theorem}
\begin {proof}
Since all involved operators are differential, it suffices to prove the
representation property locally. Hence one only has to prove
(\ref {LRepLocal}) since $\rep_\order^O$ is known to be a local
representation. But this follows directly from the above characterization
of $\DL$ and Theorem \ref{repATheo}.
\end {proof}

\begin {remark}
Though the line bundle $L$ is unique up to equivalence for
a given $B_1$, the equivalence classes
$\lambda H_{\mbox{\tiny dR}}^1(Q) / H^1_{\mbox{\tiny dR}}(Q,\Z) +
\lambda^2 H_{\mbox{\tiny dR}}^1(Q)[[\lambda]]$
of the adapted  connections $\nabla^A$ determine different
non-equivalent representations which induce `Aharonov-Bohm effects'
analogously to Sect.~\ref {ABSec}.
\end {remark}

In the case where $B_k = 0$ for $k\ge 2$, we can choose $A_k = 0$ leading
to a representation which we shall denote by $\repLo$.

Finally let us consider the case, where in addition $B$ is real.
Then $L$ carries a natural Hermitian structure $\brac{\cdot,\cdot}^L$
the equivalence classes of which are labeled by $H^1(Q,U(1))$.
Thus we can endow the space $\Gamma^\infty_{\rm cpt} (L)[[\lambda]]$
with a pre-Hilbert space structure by the $\mathbb C[[\lambda]]$-valued
Hermitian product
\begin{equation} \label{HilbertFibreMetric}
    \brac{s, t}_{L,\mu} = \int_Q \brac{s,t}^L \mu, \qquad
    s, t \in \Gamma^\infty_{\rm cpt} (L)[[\lambda]] .
\end{equation}
The following proposition now provides the adjoint of
$\repLAo (f)$ with respect to this Hermitian product.
\begin{proposition} \label{BundleAdjointProp}
  Let $B$ be in addition real and $L$ endowed with the Hermitian
  structure $\brac{\cdot,\cdot}^L$ and a formal Hermitian connection
  $\nabla^A$ obtained by a Hermitian connection $\nabla^L$ and real $A$.
  Then one has for
  $s, \tilde{s} \in \Gamma^\infty_{\rm cpt} (L)[[\lambda]]$
  and $f \in \mathcal C^\infty (T^*Q)[[\lambda]]$
\begin{equation}
    \label{wrepAdjoint}
    \brac{s, \repLAo (f)\tilde{s}}_{L,\mu}  =  \brac{\repLAo \left(
    N_{1-2\order} \cc f\right)s, \tilde{s}}_{L,\mu}.
\end{equation}
Hence $\wrepLA=\eta_{1/2}^A$ is a $^*$-representation of $\starwB$.
\end{proposition}
\begin{proof}
By $\mathbb C[[\lambda]]$-linearity and a partition of unity we may assume
that $s$ has support in some $O_j$ of a good cover of $Q$. Then the
proposition follows directly from the definition of $\brac{\cdot,\cdot}^L$,
the local expression (\ref{LRepLocal}) and Cor.~\ref{PartIntDens}.
\end{proof}
\section {GNS representation for magnetic monopoles}
\label {MagGNSSec}
In this section we consider only the Weyl ordered case that means the
product $\starwB$ for a real magnetic field
$B = \lambda B_1 + \lambda^2 B_2 + O(\lambda^3)$ satisfying the
integrality condition (IC). The complex Hermitian line bundle
$\pi_L: L \to Q$ is given as in the previous section.
Now we search for a positive linear functional such that the
induced GNS representation coincides with the
$^*$-representation $\wrepLA$.
As we shall see a slight modification will be necessary to reconstruct
$\wrepLA$ out of such a positive functional. The difficulty has its origin
in the fact that in general the line bundle $L$ does not admit a
global nowherevanishing section.
Nevertheless $L$ possesses transversal sections
$s_0 \in \Gamma^\infty (L)$. This follows from the general
transversality arguments as e.g.~in \cite[Chap. 3, Thm.~2.1]{Hirsch76}.
Now the following technical lemma should be well-known and is essential
for our constructions:
\begin{lemma} \label{TransLem}
  Let $\pi_L: L \to Q$ be a complex line bundle over $Q$ and let $\nabla$
  be a connection for $L$ and $s_0 \in \Gamma^\infty (L)$ a transversal
  section. Let $s \in \Gamma^\infty_{\rm cpt} (L)$ be a section with
  compact support. Then there exists a function
  $v \in C^\infty_{\rm cpt} (Q)$ and a vector field
  $X \in \Gamma^\infty_{\rm cpt} (TQ)$ such that
  \begin{equation} \label{TransversalSection}
    s = v s_0 + \nabla_X s_0
  \end{equation}
  and $\supp v \cup \supp X \subset \supp s$. Hereby $v$ and $X$ are not
  unique.
\end{lemma}
\begin{proof}
To avoid trivialities assume that $L$ is non-trivial and let
$Z = s_0^{-1} (\{0\}) \ne \emptyset$ be the set of zeros of $s_0$.
Choose a good cover $\{O_j\}_{j\in I}$ with corresponding local
trivializations $e_j$ of $L$. Let $x \in Z \cap O_j$ for some $j$ and
$s_0|_{O_j} = u_j e_j$ with $u_j (x) = 0$ but $du_j(x) \ne 0$
due to the transversality of $s_0$. Choose a smooth vector field $X^{(x)}$
such that $du_j (X^{(x)})|_x \ne 0$. This implies that
$\nabla_{X^{(x)}} s_0 |_x = du_j (X^{(x)}) e_j |_x \ne 0$
since $u_j(x) = 0$. By continuity there exists an open neighborhood
$U_x$ of $x$ such that $\nabla_{X^{(x)}} s_0$ does not vanish in $U_x$.
Together with $U_0 = Q \setminus Z$ the $\{U_x\}$, $x \in Z$ form an open
cover of $Q$. Hence by compactness finitely many cover $\supp s$,
say $U_0$ and $U_{x_1}, \ldots, U_{x_k}$. Let
$\{\chi_0, \chi_1, \ldots, \chi_k\}$ be a smooth partition of unity on
$U_0 \cup U_{x_1} \cup \cdots \cup U_{x_k}$ subordinate to these open
sets. Then there exist functions $v_l \in C^\infty_{\rm cpt} (Q)$ for
$l =0, \ldots, k$ such that
$ \chi_0 s = v_0 s_0$, $\chi_l s = v_l \nabla_{X^{(x_l)}} s_0$
and  $\supp v_l \subset \supp s \cap U_l$. Then we have globally
$s = v_0 s_0 + \nabla_X s_0$ with
$X = v_1 X^{(x_1)} + \cdots + v_k X^{(x_k)}$. This proves the claim.
\end {proof}

Due to this lemma we can reconstruct any section of $L$ with compact
support out of a fixed transversal section $s_0$ if we allow for covariant
differentiation. Now the idea is that with $\wrepLA$ we have indeed the
possibility of covariant differentiation but this automatically raises the
degree in $\lambda$. Hence we have to allow finitely many negative
powers of $\lambda$ to get covariant differentiation of $s_0$ in the
lowest order in $\lambda$. In other words we have to replace formal power
series by formal Laurent series everywhere and thus obtain vector spaces
over the field of formal Laurent series $\LS{\mathbb C}$ instead of
modules over $\mathbb C[[\lambda]]$. This does not cause any problems as
can be seen directly from all constructions above or the general
arguments given in \cite [App.~A]{BNW97b}.

Thus having fixed a transversal section $s_0 \in \Gamma^\infty (L)$ we
define the $\LS{\mathbb C}$-linear functional
$\omega_{s_0}: \LS{\mathcal C^\infty_Q (T^*Q)} \to \LS{\mathbb C}$ by
\begin{equation} \label{omegaTransDef}
    \omega_{s_0} (f) := \int_Q  \brac{ s_0, \wrepLA (f) s_0}^L \mu
\end{equation}
which is clearly well-defined for $f \in \LS{\mathcal C^\infty_Q (T^*Q)}$.
The positivity of $\omega_{s_0}$ is established easily.
\begin{lemma} \label{TransOmegaPosLem}
The $\LS{\mathbb C}$-linear functional $\omega_{s_0}$ is positive and the
Gel'fand ideal $\mathcal J_{s_0}$ is given by all $f$
fulfilling $\wrepLA (f) s_0 = 0$. Furthermore one has for
$f, g \in \LS{\mathcal C^\infty_Q (T^*Q)}$
\begin{equation}\label {TransOmegaFactor}
    \omega_{s_0} (\cc f \starwB g) =
    \int_Q \brac{ \wrepLA (f) s_0, \wrepLA (g) s_0}^L \mu .
\end{equation}
\end {lemma}
\begin {proof}
Clearly it suffices to prove (\ref{TransOmegaFactor}). Due to the
obvious $\LS{\mathbb C}$-linearity we may consider
$f, g \in \mathcal C^\infty_Q (T^*Q)$ without $\lambda$-powers.
In order to apply Prop.~\ref{BundleAdjointProp}, where we need
sections with compact support, we choose a smooth function
$\phi: Q \to [0,1]$ with compact support such that $\phi$ is identical
to $1$ on  $\supp (\iota^*f) \cup \supp(\iota^*g)$.
Due to the representation property and (\ref{wrepAdjoint})
\[
\begin{split}
   \omega_{s_0} (\cc f \starwB g)
   & = \int_Q \brac{ s_0, \wrepLA (\cc f \starwB g) s_0}^L \mu
     = \int_Q \brac{ \phi s_0, \wrepLA (\cc f \starwB g) \phi s_0} \mu \\
   & = \int_Q \brac{ \wrepLA (f) \phi s_0, \wrepLA (g) \phi s_0}^L \mu
     = \int_Q \brac{ \wrepLA (f) s_0, \wrepLA (g) s_0} \mu
\end{split}
\]
which proves (\ref{TransOmegaFactor}), thus the lemma.
\end {proof}

The following theorem shows that $\wrepLA$ can indeed be recovered as
GNS representation induced by $\omega_{s_0}$. Note that the
transversality as well as the usage of formal Laurent series is crucial.
\begin{theorem} \label {FormalTransTheo}
  Let $s_0 \in \Gamma^\infty (L)$ be a transversal section and
  $\omega_{s_0}$ as in (\ref{omegaTransDef}). Then the GNS representation
  space $\mathfrak H_{s_0} =
  \LS{\mathcal C^\infty_Q (T^*Q)}\big/\mathcal J_{s_0}$
  is canonically isometric to $\LS{\Gamma^\infty_{\rm cpt} (L)}$ via
  \begin{equation} \label{BundleRepIso}
    \mathfrak H_{s_0} \ni \psi_f \mapsto \wrepLA (f) s_0 \in
    \LS{\Gamma^\infty_{\rm cpt} (L)},
  \end{equation}
  where $\LS{\Gamma^\infty_{\rm cpt} (L)}$ is equipped with the Hermitian
  product (\ref{HilbertFibreMetric}). Moreover, the GNS representation on
  $\mathfrak H_{s_0}$ is unitarily equivalent to $\wrepLA$ under the
  unitary map (\ref {BundleRepIso}) and extends to a representation of
  $\LS{\mathcal C^\infty (T^*Q)}$.
\end{theorem}
\begin{proof}
First notice that $\psi_f \mapsto \wrep^L (f) s_0$ is well-defined and
isometric since due to Lemma \ref{TransOmegaPosLem}
\[
    \brac{\psi_f,\psi_g} = \omega_{s_0} (\cc f \starwB g)
    = \int_Q \brac{ \wrepLA (f)s_0, \wrepLA (g) s_0}^L \mu
\]
and thus (\ref{BundleRepIso}) is injective. To prove surjectivity we
apply Lemma \ref{TransLem}: let first
$s \in \LS{\Gamma^\infty_{\rm cpt} (L)}$ be a section without $\lambda$
powers and extend the arguments by $\LS{\mathbb C}$-linearity afterwards.
Then due to the lemma there exists a function
$v_0 \in \mathcal C^\infty_{\rm cpt} (Q)$ and a vector field
$X_0 \in \Gamma^\infty_{\rm cpt} (TQ)$ such that
$s = v_0 s_0 + \nabla_{X_0} s_0$. Due to the explicit formula
(\ref{LwrepDef}) for the standard order representation $\srepLA$
we find that $\srepLA (\pi^* v_0 + \frac{\im}{\lambda} J(X_0)) s_0$
coincides with $s$ up to higher orders in $\lambda$. By iterating this
procedure we find a function
$v \in \mathcal C^\infty_{\rm cpt} (Q)[[\lambda]]$ and a vector
field $X \in \Gamma^\infty_{\rm cpt} (TQ)[[\lambda]]$ such that
$\srepLA (\pi^* v + \frac{\im}{\lambda} J(X)) s_0 = s$. This proves
the surjectivity since $\srepLA$ can be replaced by $\wrepLA$
incorporating the bijection $N_{1/2} (\alpha_\mu)$.
The other statements of the theorem follow immediatly.
\end{proof}

\begin{remark}
In the context of formal GNS constructions the usage of
the field of formal Laurent series was extensively discussed in
\cite{BW98a}. Furthermore it may be desirable to use even an
algebraically closed extension field for these constructions. Hereby
the field of \emph{formal completed Newton-Puiseux series} $\CNP{\mathbb C}$
is the natural candidate. We mention that anything can also be done in this
framework due to the `extension theorems' in \cite [App.~A] {BNW97b}.
\end{remark}

\section{Symbol calculus and star products}
\label{SymbSec}
In the previous sections we have constructed representations $\repAo$
of the deformed product $\starBo$ on formal sections of line bundles $L$
over $Q$. We will show in the following that the formal
representation $\repAo$ can be interpreted as the asymptotic expansion
of a certain pseudodifferential operator calculus on $Q$ with values in $L$.
In other words, every $\repAo$ comes from a particular global symbol
calculus on $Q$. So finally we obtain an analytical interpretation
for the formal representations $\repAo$.

But before we go into the details of the construction we have to set up
some notation and will recall some results from the theory of
pseudodifferential operators on  manifolds. For more information on this
we refer the reader to {\sc H\"ormander} \cite{Hor:ALPDOIII}.

Under $\sym^m (Q)$, $m \in \R$ we understand the symbols on $T^*Q$ of
H\"ormander type, i.e.~$\sym^m (Q)$ consists of all smooth functions
$a \in \mathcal C^{\infty} (T^*Q)$ such that uniformly on compact subsets
$K \subset U $ of any coordinate patch $U \subset Q$
\begin{equation} \label{SymbEst}
  \left| \partial^{\alpha}_{q^\cdot} \, \partial^{\beta}_{p_\cdot} \,
  a (\zeta) \right|
  \, \leq C_{K,\alpha,\beta} \, \brac{\zeta}^{m-|\beta|},
  \qquad  \zeta \in T^*K, \quad \alpha,\beta \in \N^n.
\end{equation}
Hereby $\brac{\zeta} = ( 1 +  \|p(\zeta)\|^2)^{1 /2} $ and
$ C_{K,\alpha,\beta} >0 $ and we use the usual multi-index notation
for $\alpha, \beta$. Note that this definition is independent of the used
chart and provides a global characterization of symbols. As usual we set
$\sym^{\infty} (Q) = \cup_{m \in \R} \, \sym^m (Q)$
and $\sym^{- \infty} (Q) = \cap_{m \in \R} \, \sym^m (Q)$.
In case $E\rightarrow Q$ is a smooth vector bundle one defines
$\sym^m (Q;E)$ as the space of all sections
$a \in \Gamma^\infty (\pi^* (E))$ such that for every
section $e \in \Gamma^\infty (E^*)$ of the dual bundle of $E$ the relation
$\brac{e,a} \in \sym^m (Q)$ holds. The spaces ${\mathcal P} (Q; E)$
and  ${\mathcal P}^k (Q; E)$ of fiberwise homogeneous resp.~fiberwise
homogeneous smooth sections of degree $k$ of $\pi^* (E)$ are then
natural subspaces of $\sym^\infty (Q; E)$ resp. of $\sym^k (Q; E)$

For two smooth vector bundles $E\rightarrow Q$ and $F \rightarrow Q$ the
space $\Psi^m (Q;E,F)$ of pseudodifferential operators from $E$ to $F$ of
order $m$ consists of all continuous mappings
$A: \Gamma^\infty_{\rm cpt} (Q,E) \rightarrow \Gamma^\infty_{\rm cpt} (Q,F)'$
such that the Schwartz kernel
$K^{\! A} \in \Gamma^\infty_{\rm cpt} (Q \times Q , F \boxtimes E^* )'$
is a conormal distribution  of order $m$  with respect to the diagonal of
$Q \times Q$ (cf.~\cite[Def.~18.2.6.]{Hor:ALPDOIII}).
A global symbol calculus now provides an (almost) bijective correspondence
between symbols and  pseudodifferential operators over $Q$.
The symbol calculus we introduce in the following generalizes the
well-known ones of {\sc Widom} \cite{Wid:CSCPO} and
{\sc Safarov} \cite{Saf:POLC}, so we can be brief with proofs
(cf.~also \cite{Pfl:NSRM,Vor:unpublished}).

By a {\it cut-off function} for the connection $\nabla$ on $Q$ we mean
a smooth function $\chi : TQ \rightarrow [0,1]$ having support
in an open neighboorhood $O \subset TQ$ of the zero section such that
$(\pi,\exp)|_O$ is a diffeomorphism onto an open neighborhood $\widetilde{O}$
of the diagonal of $Q\times Q$. Additionally let us suppose that
$\supp \chi \cap T_xQ$ is compact for every $x \in Q$. After shrinking $O$ we
can assume that for every $\order \in [0,1]$ the map
$\chi_\order : TQ \rightarrow [0,1]$ with
\begin{equation}\label{PsiOrdDef}
  \chi_\order (v_x) = \left\{
  \begin{array}{lll}
    \chi (-\order T_{\order v_x} \exp_x v_x)
    \chi ( (1-\order ) T_{\order v_x} \exp_x v_x) & \mbox{ for } &
    v_x \in T_xQ \cap  O \\
    0 &\mbox{ else }
  \end{array} \right.
\end{equation}
is smooth and itself a cut-off function.
The volume densities $\mu_x$ and $\mu$ are now related by a smooth function
$\rho : Q\times Q \supseteq \widetilde O \to \mathbb R^+$ satisfying
\begin{equation}
    \big(\exp_x^{-1}\big)^* (\mu_x) = \rho(x, \cdot) \mu.
\end{equation}

Next assume to be given vector bundles $E$ and $F$ over $Q$, a Hermitian
metric $\brac{\cdot,\cdot}^F$ on $F$ and covariant derivatives
$\nabla^E$ and $\nabla^F$, where the latter one is supposed to preserve
the Hermitian  metric.
Then, if $x, y \in Q$ are sufficiently close to each other, or in
other words if $(x, y) \in \widetilde O$, we will denote by
$\TAU {E} {y} {x}: E_y = \pi^{-1}(y) \to E_x$
the parallel translation in $E$ along the unique geodesic with respect
to $\nabla$ connecting $y$ with $x$ in $\exp_y(T_y Q \cap O)$.
Moreover, $\ormid {x} {y}= \exp_x(\order \exp_x^{-1}(y))$ will denote
the \emph{$\order$-midpoint} of this geodesic.
Now define for every section $e \in \Gamma^\infty(E)$ the lift
$e^\chi \in \Gamma^\infty (\pi^*_{{}_{TQ}}(E))$ by
\begin{equation} \label {chiLiftDef}
  e^\chi(v_x)=\left\{
  \begin{array} {lll}
    \chi(v_x) \TAU {E} {\exp_x v_x} {x}
    (e(\exp_x v_x)) & \mbox{ for all } & v_x \in T_xQ \cap O, \\
    0 & \mbox{ else, }
  \end{array} \right.
\end{equation}
where $x$ runs through $Q$.
So finally we have all the ingredients to formulate, how
pseudodifferential operators are associated to symbols.
\begin{theorem} \label {ThmPseudoOpDef}
    Let $\order \in [0,1]$ be an ordering parameter and
    $\hbar \in \R^+$ a deformation parameter.
    Define for every symbol $a \in \sym^m (Q; \Hom (E,F))$ and all sections
    $f \in \Gamma^\infty (F)$ and $e \in \Gamma^\infty (E)$ the
    kernel function $W_\order (a)(f, e) \in \mathcal C^\infty(Q)$ by
    \begin{eqnarray}
       \big( W_\order (a)(f, e)\big)(x)=
       \frac{1}{(2\pi \hbar)^n} \int_{T_x^* Q} \int_{T_x Q}
       \eu^{-\frac{\im}{\hbar} \brac{\zeta_x,v_x }}
       \brac{ f^\chi( -\order v_x), a (\zeta_x)
       (e^\chi ((1-\order)v_x))}^F \mu_x (v_x)
       {\mu_x}^{\!\!\! *}(\zeta_x),\nonumber
    \end{eqnarray}
   where $\chi $ is a cut-off function as introduced above.
   Then there exists a unique properly supported pseudodifferential operator
   $\Opo (a) \in \Psi^\infty (Q;E,F)$ of order $m$ fulfilling the relation
   \begin{equation}
     \brac{f, \Opo (a) \,e  }
     = \int_Q \, \brac{f(x),(\Opo (a)\, e)(x)}^F \,
        \mu (x)
     = \int_Q \, \big( W_\order (a) (f,e) \big) (x) \, \mu (x)
   \end{equation}
   for all $f$ and $e$ with compact support.  One calls
   $\Opo (a)$ the $\kappa$-ordered quantization of $a$.

   If now $\chi'$ is another cut-off function the difference
   $\Op_{\hbar,\order}^{[\chi]} (a) - \Op_{\hbar,\order}^{[\chi']} (a)$
   is a smoothing operator, hence lies in $\Psi^{-\infty} (Q;E,F)$.
   In case $a$ is polynomial in the momenta, i.e.~if
   $a \in {\mathcal P} (Q;\Hom(E,F))$,
   $\Opo (a)$ is even a differential operator between $E$ and $F$ and
   independent of the choice of $\chi$.
\end{theorem}
\begin{remark}
In a few cases we need to specify the dependence of the operator calculus
on the choice of connections $\nabla$ on the considered vector bundles
or the cut-off function $\chi$. When necessary, we will denote this
dependence by $\Op_{\hbar,\order}^{[\nabla,\chi]}$.
\end{remark}
\begin{proof}
   For the case, where $E$ and $F$ are trivial line bundles over $Q$, the
   claim has been shown in {\sc Pflaum} \cite[Thm.~2.1]{Pfl98c}.
   By local triviality of $E$ and $F$ the general statement then follows
   easily.
\end{proof}
It turns out that the quantization map
$\Opo: \sym^\infty (Q;\Hom (E,F)) \rightarrow \Psi^\infty (Q;E,F)$
possesses a filtration preserving pseudo-inverse
$\sigma_\order : \Psi^\infty (Q;E,F)\rightarrow \sym^\infty (Q;\Hom (E,F))$
called a {\it symbol map}.
More precisely this means that the induced maps
$ \cc \Opo : \sym^\infty / \sym^{-\infty} (Q;\Hom (E,F)) \rightarrow
\Psi^\infty / \Psi^{-\infty} (Q;E,F)  $ and
$\cc \sigma_\order :  \Psi^{\infty} / \Psi^{-\infty} (Q;E,F) \rightarrow
\sym^{\infty} / \sym^{-\infty} (Q;\Hom(E,F)) $
are inverse to each other. Let us give an explicit representation for
the symbol map $\sigma_\order$.
So consider a pseudodifferential operator $A \in \Psi^m (Q; E,F)$
and denote by $K^{\! A}_\order$  the $\pi_{TQ}^* (\Hom (E,F))$-valued
distribution defined by
\begin{equation}
   K^A_\order (v_x)=\chi(v_x) \TAU {F} {\exp_x (-\order v_x)} {x}
   (K^A \rho^{-1})(\exp_x (-\order v_x),\exp_x ((1-\order) v_x))
   \TAU {E} {x} {\exp_x((1-\order)v_x)}
\end{equation}
for all $v_x \in T_xQ \cap  O$ and set to zero outside $O$.
Then the symbol $\sigma_\order (A)$ is given by the oscillatory integral
\begin{equation}
    \sigma_\order (A) (\zeta_x) = \int_{T_x Q} \eu^{\frac{\im}{\hbar}
    \brac{ \zeta_x, v_x }} K^A_\order (v_x) \, \mu_x(v_x).
\end{equation}
In the following we will provide some useful integral representations
for the operator $\Opo (a)$. The proof is done by performing straightforward
but somewhat lenghty transformations of the involved integrals and
applying Gau{\ss} Lemma,
i.e. $T_{v_x} \exp_x v_x = \TAU { }{x}{\exp_x (v_x)} v_x$,
several times.
\begin{proposition}\label{IntRepProp}
For every symbol $a \in \sym^m( Q; \Hom (E, F))$, $m\in \mathbb R$ and every section
$e\in \Gamma^\infty_{\rm cpt} (E)$, $\Opo(a)e$ can be written in the form
\begin{equation}
\begin{split}
\lefteqn{(\Opo(a)e)(x)} \\
  & = \frac{1}{(2\pi \hbar)^n} \int_{T_x^*Q} \int_{T_xQ}
  \eu^{-\frac{\im}{\hbar} \brac{\zeta_x, v_x}} \, \TAU {F}
  {\exp_x(\order v_x)} {x} a\left(
  (T_{\exp_x (\order v_x)}\exp_x^{-1})^* \zeta_x \right)
  \TAU {E} {\exp_x v_x} {\exp_x (\order v_x)}
  e(\exp_x(v_x))
  \\  &
  \hspace{15mm} \rho(\exp_x (\order v_x), x)
  \, \chi_\order (v_x) \, \mu_x (v_x) \, \mu_x^*(\zeta_x).
\end{split}
\end{equation}
Moreover, the Schwartz kernel $K^{\!A} $ of $A =\Opo(a)$ is given by
\begin{equation}
\begin{split}
  K^{\! A} (x, y) &= \frac{1}{(2\pi \hbar)^n} \int_{T_x^*Q}
  \eu^{-\frac{\im}{\hbar} \brac{\zeta_x, \exp_x^{-1}(y)}}\,
  \TAU {F} {\ormid {x} {y}} {x} a\left( (T_{\ormid {x} {y}}
  \exp_x^{-1})^*\zeta_x\right) \TAU {E} {y} {\ormid {x} {y}}\\ \
  & \hspace{25mm} \rho(\ormid{x}{y} , x) \, \rho( x, y) \,
  \widetilde{\chi_\order}(x, y) \, \mu_x^*(\zeta_x),
\end{split}
\end{equation}
where $\widetilde{\chi_\order}$ is the function equal to
$\chi_\order \circ (\pi_{_{TQ}},\exp) $ over $\widetilde{O}$ and which is
extended by zero on the complement of $\widetilde{O}$ in $Q\times Q$.
\end{proposition}

Next we will calculate more explicitly the differential operator
$\Opo ({\mathsf p})$ respectively its kernel function
$W_\order ({\mathsf p}) (f,e) $ for a vector valued symbol
${\mathsf p} \in {\mathcal P}^k (Q;\Hom(E,F))$
polynomial in the momenta of degree $k$.
In local coordinates one can write ${\mathsf p}$ in the form
\[
  {\mathsf p}
  = \frac{1}{k!} \sum_{j_1, \cdots , j_k} \, p_{j_1} \cdots p_{j_k} \,
  T^{j_1 \cdots j_k} \circ \pi ,
\]
where $T^{j_1 \cdots j_k} \in \Gamma^\infty (\Hom(E,F))$.
Now denote by  $v^j$ the local coordinates on $T_xQ$ induced by $q$, and
by $V_x : \R^n \rightarrow T_xQ$ the function
$(v^1,\ldots,v^n)\mapsto
v^1 \partial_{q^1}\big|_{x} + \cdots + v^n \partial_{q^n}\big|_{x}$,
where $x \in Q$ lies in the domain of $q$.
Then we compute by integration by parts
\[
\begin{split}
\lefteqn{\left( W_\order({\mathsf p}) (f,e)\right) (x)  = }\\
  & =
  \frac{(-\im\hbar)^k}{k!} \sum_{j_1,\ldots,j_k}
  \partial_{v^{j_1}} \cdots \partial_{v^{j_k}}
  \Big|_{v=0} \brac{ f^\chi (- \order V_x(v)) ,
  T^{j_1 \cdots j_k} (x) \left( e^\chi((1-\order )V_x(v))\right)}_x^F
  \\
  & =
  \frac{(-\im \hbar)^k}{k!} \sum_{j_1,\ldots,j_k}
  \sum_{l=0}^k \, \left( k \atop l \right) (-\order)^l (1-\order)^{k-l} \\
  & \hspace{1.5em} \brac{
  \partial_{v^{j_1}} \!\cdots \partial_{v^{j_l}} \Big|_{v = 0}
  \TAU {F} {\exp_x V_x} {x}( f(\exp_x V_x )) , T^{j_1 \cdots j_k} (x)
  \left( \partial_{v^{j_{l+1}}} \!\cdots \partial_{v^{j_k}} \Big|_{v =0}
  \TAU {E} {\exp_x V_x} {x} ( e (\exp_x V_x )) \right) }_x^F .
\end{split}
\]
Using the symmetrized covariant derivate $\mathsf D$ one checks that
\begin{equation}
  \partial_{v^{j_1}} \!\cdots \partial_{v^{j_l}} \Big|_{v=0}
  \TAU {F} {\exp_x V_x} {x}( s (\exp_x V_x )) = \frac{1}{l!}
  \brac{ \partial_{q^{j_1}} \otimes \cdots \otimes
  \partial_{q^{j_l}} , {\mathsf D}^l s } (x).
\end{equation}
Inserting this relation one finally obtains
\begin{equation} \label{WignerKRep2}
\begin{split}
  \left( W_\order({\mathsf p}) (f,e)\right)
  = & \frac{(-\im \hbar)^k}{k!} \sum_{j_1,\ldots,j_k} \sum_{l=0}^k \,
  \left( k \atop l \right) (-\order)^l (1-\order)^{k-l} \frac{1}{l!}
  \frac{1}{(k-l)!} \\
  &
  \brac{ \brac{ \partial_{q^{j_1}} \otimes \cdots \otimes \partial_{q^{j_l}} ,
  {\mathsf D}^l f } , T^{j_1 \cdots j_l} \left(
  \brac{ \partial_{q^{j_{l+1}}} \otimes \cdots \otimes \partial_{q^{j_k}},
  {\mathsf D}^{k-l} e } \right) }^F .
\end{split}
\end{equation}
Now, let us abbreviate for every
$T \in \Gamma^\infty (\bigvee^k TQ \otimes \Hom(E,F))$ and every $l \leq k$
by
$T_l: \Gamma^\infty (E) \rightarrow \Gamma^\infty ( \bigvee^l TQ \otimes F)$
the differential operator which locally is given by
\[
  T_l (e) = \frac{1}{(k-l)!} \sum_{j_{l+1}, \cdots , j_k} \,
  \brac{ dq^{j_{l+1}} \otimes \cdots \otimes dq^{j_k} , T
  \left( \brac{\partial_{q^{j_{l+1}}} \otimes \cdots \otimes
  \partial_{q^{j_k}}, {\mathsf D}^{k-l} e }\right)}.
\]
Then we have the following result.
\begin{theorem}
  Let ${\mathsf p} \in {\mathcal P}^k (Q ; \Hom (E,F))$ and
  $e \in \Gamma^\infty_{\rm cpt} (E)$.
  Then
  \begin{equation}
    \Opo ({\mathsf p} ) e = \sum_{l=0}^k \, \frac{(-\im \hbar)^k}{l!(k-l)!}
    \order^l  (1-\order)^l \, \left(
    \div_\mu^F\right)^l \left( T_l (e) \right) ,
  \end{equation}
  where $T \in \Gamma^\infty (\bigvee^k TQ \otimes \Hom(E,F))$ is the unique
  tensor field fulfilling ${\mathsf p} = J (T)$.
\end{theorem}
\begin{proof}
  The claim follows from Eq.~(\ref{WignerKRep2}) by Lemma
  \ref{LemAdjointNabla} and the definition of the operator
  $J$.
\end{proof}
\begin{corollary}\label{CorOpChang}
  For every ${\mathsf p} \in {\mathcal P} (Q ; \Hom (E,F))$ one has
  \begin{equation} \label {OpoPoly}
    \Opo ({\mathsf p}) = \Op_0 \left( J \exp\left(- \im \order \hbar \,
    \div_\mu^{\Hom(E,F)} \right) J^{-1} ({\mathsf p}) \right)
    = \Op_0 \left(\exp\left(-\im\order\hbar\Delta_\mu^{\Hom (E, F)} \right)
      {\mathsf p}\right).
  \end{equation}
\end{corollary}

Now let us come back to our original goal, namely to give an operator
theoretical interpretation of the representations $\repAo$
on the formal sections of a line bundle $L$ over $Q$.
So let $E = F = L$ be a Hermitian line bundle
with metric connection $\nabla^L$ associated to a closed formal real
two form $ B = \lambda B_1$ satisfying (IC).
Comparing formula (\ref {OpoPoly}) with Eq.~(\ref {LwrepDef}) we obtain
the following second corollary.
\begin {corollary}
For every integral real formal two-form $B = \lambda B_1$ and every
corresponding line bundle $L$ with Chern class $[\frac{1}{2\pi} B_1]$
one has the relation
\begin {equation} \label {OpoRepLo}
    \Opo (\mathsf p) = \left. \repLo (\mathsf p)\right|_{\lambda = \hbar},
    \qquad
    \mathsf p \in \mathcal P (Q; \Hom (L, L)) = \mathcal P (Q),
\end {equation}
where we have substituted the formal variable $\lambda$ in
$\repLo (\mathsf p)$ by the real value $\hbar \in \mathbb R^+$.
\end {corollary}
\begin {remark}
By formula (\ref {LwrepDef}), $\repLo (\mathsf p)$ is a polynomial
in $\lambda$, hence insertion $\hbar$ for $\lambda$ is well-defined indeed.
\end {remark}

By the last corollary one can interpret the pseudodifferential calculus
$\Opo$ as an operator realization for the formal representation $\repLo$.
Let us now briefly recall that for two symbols
$a \in \sym^m (Q; \Hom (L, L))$ and
$b \in \sym^{\tilde m} (Q; \Hom (L, L))$ the pseudodifferential operator
$\Opo (a) \Opo (b)$ lies in $\Psi^{m+\tilde m} (Q; L, L)$. Additionally,
over every coordinate patch $U$ of $Q$ trivializing $L$ one even has
an asymptotic expansion of the symbol $\sigma_0 (\Op_0 (a) \Op_0 (b))$ of
the form
\begin {equation} \label {SymbolAsympt}
    \left.\sigma_0 (\Op_0 (a) \Op_0 (b))\right|_{T^*U}
    \sim
    \sum_{k=0}^\infty \hbar^k\sum_{{|\alpha|, |\beta|+|\gamma|\leq
    k \atop
    |\alpha|+|\beta|-|\delta|= k}}
    \pi^* (r^\delta_{k,\alpha\beta\gamma})
    p_\delta \partial^\alpha_{p_{\cdot}} a
    \partial^{\beta}_{p_{\cdot}} \partial^\gamma_{q^{\cdot}} b ,
\end {equation}
where $r^\delta_{k,\alpha\beta\gamma} \in \mathcal C^\infty (U)$ and $\alpha,
\beta, \gamma, \delta \in \mathbb N^n$ denote multi-indices.
(cf.~\cite[Sect.~5]{Wid:CSCPO} and \cite[Thm.~6.4]{Pfl:NSRM}).
By the corollary and relation (\ref{SymbolAsympt}) we can now conclude
that for all symbols $a \in \sym^m (Q)$ and
$b \in \sym^{\tilde m} (Q)$ there is an asymptotic expansion of the
operator product $\Opo (a) \Opo (b)$ of the form
\begin{equation} \label{OpPropAsympt}
  \Opo (a) \Opo (b) \sim
  \Opo (ab) + \sum_{k=1}^\infty \hbar^k \Opo (C_k (a, b)),
\end{equation}
where the $C_k$ denote the $k$-th bidifferential operator of the star product
$\starBo$. Moreover Eq.~(\ref {OpoRepLo}) shows that $\repLo$ is
induced by $\Opo$.

In the following we want to generalize this result to the case of the star
product $\starBo$ when $B = \lambda B_1 + \lambda^2 B_2 + O(\lambda^3)$ is
real formal two-form fulfilling (IC). Now, choose  a formal vector potential
$A = \lambda^2 A_2 + O(\lambda^3)$  like in
Section~\ref{MagMonoSec} such that $dA = B - \lambda B_1$.
Next we choose a smooth map
$\mathbb R^+ \to \Omega^1 (Q), \hbar \mapsto A_\hbar$ having Taylor
expansion $A$ around $\hbar = 0$. Then let
$B_\hbar = \hbar B_1 + dA_\hbar$, and
$\nabla^L_\hbar = \nabla^L + \frac{\im}{\hbar} A_\hbar$ and check that the
curvature of $\nabla^L_\hbar$ is just $\im B_\hbar$.
Finally we associate to every symbol $a \in \sym^\infty (Q)$ and all
$\hbar \in \mathbb R^+$, $\order \in [0,1]$ the pseudodifferential
operator $\Op_{\hbar,\order}^A (a) =
\Op^{[\nabla^L_\hbar,\chi]}_{\hbar,\order}(a)$.
Then we arrive at the following main result.
\begin{theorem}\label{AsymExpTheo}
For all real formal series
$B = \lambda B_1 + \lambda^2 B_2 + O(\lambda^3)$ and
$A = \lambda^2 A_2 + O(\lambda^3)$ of two- resp.~one-forms such that
$B$ is integral and $B = \lambda B_1 + dA$ the
operator map $\Op^A_{\hbar,\order}$ comprises an operator representation
of $\starBo$ on the Hermitian line bundle $L$ with Chern class
$\frac{1}{2\pi} [B_1]$ and metric connection $\nabla^L$. Moreover the
formal representation $\repLAo$ of $\starBo$ is induced by the operator
calculus $\Op^A_{\hbar,\order}$. In case $A=0$ one has
$\Op^A_{\hbar,\order} = \Op_{\hbar,\order}$.
\end {theorem}
\begin{proof}
The relation $\Op^A_{\hbar,\order} = \Opo$ is obvious in case $A=0$, hence
the claim follows for $A=0$. Now, fix $\hbar$ and let
$\hat{\hbar} \in \mathbb R^+$ be arbitrary. Then again by the above
considerations the formal representation $\tilde \eta_\order$ of
$\star_\order^{\frac{\lambda}{\hbar} B_\hbar}$ defined by the Hermitian
connection $\nabla^L_\hbar$ and the formal one-form $\tilde A = 0$ is
induced by the operator calculus
$\Op_{\hat\hbar,\order}^{[\nabla^L_\hbar]}$, or formally
\[
    \Op_{\hat\hbar, \order}^{[\nabla^L_\hbar]} (a) =
    \sum_{l=0}^\infty {\hat\hbar}^l \tilde \eta_{\order,l} (a)
    \bmod \hat\hbar^\infty,
\]
where $\tilde \eta_{\order,l}$ depends on $\hbar$.
From Eq.~(\ref {wrepAFormel}) it follows that one has
$\sum_{l=0}^\infty \hbar^l \tilde \eta_{\order,l}(a)
=\sum_{l=0}^\infty \hbar^l \eta^A_{\order,l} (a)$ for $\hbar = \hat\hbar$.
Consequently
\[
    \Op^A_{\hbar,\order} (a) = \sum_{l=0}^\infty
    \hbar^l \eta^A_{\order,l} (a)
    \bmod \hbar^\infty
\]
which entails the claim.
\end{proof}

\section {GNS construction for pseudodifferential operators}
\label {GNSSymbSec}
In this section we show how the symbol and pseudodifferential operator
calculus for the Weyl ordered type, i.e.~for $\order = 1/2$ can be
understood as a particular GNS construction analogously to the formal
case as considered in Theorem \ref {FormalTransTheo}.

To achieve this we will first introduce an integral representation of the
formal operator $N_\order$. Hereby we can assume without restriction that
the cut-off function $\chi$ is the square of another cut-off function
$\sqrt{\chi}$, and likewise for $\chi_\order$. Moreover, we denote by
$\csym^m (Q; \Hom (E,F))$ those symbols $a$ in $\sym^m (Q; \Hom(E,F))$
such that $\cc{\pi(\supp a)} \subset Q$ is compact.
Then the following proposition determines the relation between the
different orderings in the various operator calculi.
\begin {proposition} \label {NopoProp}
\begin {enumerate}
\item The operator
      $\Nopo: \sym^\infty (Q; \Hom (E, F)) \to
      \sym^\infty (Q; \Hom (E, F))$ defined by
      \begin {equation} \label {NopoDef}
          \begin {split}
          (\Nopo a) (\zeta_x') = &
          \frac{1}{(2\pi\hbar)^n} \int_{T^*_x Q} \int_{T_xQ}
          \eu^{-\frac{\im}{\hbar} \brac{\zeta_x, v_x}}
          \sqrt{\chi_\order} (v_x)
          \rho (\exp_x (\order v_x), x) \\
          &
          \TAU{F}{\exp_x (\order v_x)}{x}
          a \left(
          (T_{\exp_x (\order v_x)}\exp_x^{-1})^* (\zeta_x'+\zeta_x)
          \right)
          \TAU{E}{x}{\exp_x (\order v_x)}
          \mu_x (v_x) {\mu_x}^{\!\!\! *}(\zeta_x)
          \end {split}
      \end {equation}
      for $a \in \sym^\infty (Q; \Hom (E, F))$ maps
      $\csym^m (Q;\Hom(E,F))$ to $\csym^m (Q; \Hom(E,F))$ and
      fulfills the equation
      \begin {equation} \label {OpoOpNopo}
          \Op_\order^{[\chi]} (a) = \Op_0^{[\sqrt{\chi_\order}]} (\Nopo a) .
      \end {equation}
\item There exists a cut-off function $\psi_\order: TQ \to [0,1]$ such that
      \begin {equation} \label {OpOpoCCNopo}
          \Op_0^{[\psi_\order]} (a) = \Op_0^{[\sqrt{\chi_\order}]}
          (\Nopo \cc{\Nopo} a) = \Op_\order^{[\chi]} (\cc{\Nopo} a)
      \end {equation}
      for every $a \in \sym^m (Q; \Hom (E, F))$. Hereby $\cc \Nopo$ is
      defined as $\Nopo$ by Eq.~(\ref{NopoDef}) only with
      complex conjugated phase factor
      $\eu^{\frac{\im}{\hbar}\brac{\zeta_x,v_x}}$.
      In particular $\Nopo$ and $\cc\Nopo$ are quasi-inverse to each
      other, i.e.~$\Nopo\cc\Nopo - \id$ and $\cc\Nopo\Nopo - \id$ are
      smoothing.
\item For every $\mathsf p \in {\mathcal P} (Q ; \Hom (E,F))$ and all
      $\order \in [0,1]$ one has
      \begin {equation}
          \Nopo \mathsf p =
          \exp \left(- \im \order \hbar
          \Delta_\mu^{\Hom (E,F)} \right) \mathsf p
          \quad
          \textrm{ and }
          \quad
          \cc \Nopo \mathsf p =
          \exp \left(\im \order \hbar
          \Delta_\mu^{\Hom (E,F)} \right) \mathsf p.
      \end {equation}
      In general $\Nopo$ and $\cc \Nopo$ possess Taylor expansions
      which are given explicitly by
      \begin {equation}
          \Nopo a =
          \exp \left(- \im \order \hbar
          \Delta_\mu^{\Hom (E,F)} \right) a \bmod \hbar^\infty
          \quad \textrm { and } \quad
          \cc{\Nopo} a =
          \exp \left( \im \order \hbar \Delta_\mu^{\Hom (E,F)} \right) a
          \bmod \hbar^\infty.
      \end {equation}
\end {enumerate}
\end {proposition}
\begin {proof}
For the proof of Eq.~(\ref {OpoOpNopo}) first check the claim on
$\sym^{-\infty}(Q; \Hom(E,F))$ by straightforward computation and extend
it to $\sym^\infty (Q; \Hom(E,F))$ by a continuity argument. Likewise one
shows part ii.). The fact that $\sym^m_{\rm cpt} (Q;\Hom (E,F))$ is mapped
into itself under $\Nopo$ follows from the fact that $\cc{\pi (\supp a )}$
and $\supp (\chi_\order|_{T_x Q})$ are compact and the integral representation
(\ref{NopoDef}). Corollary \ref {CorOpChang} and part i.) entail part
iii.) by first checking the claim on polynomial symbols
$\mathsf p \in \mathcal P (Q; \Hom(E, F))$.
\end {proof}

From now on we specialize to the case of Weyl ordered calculus and
$E=F=L$, where $L$ is a Hermitian line bundle with Hermitian connection.
The following two lemmas are crucial for the GNS construction we have
in mind.
\begin {lemma}
Let $s_0 \in \Gamma^\infty (L)$ be a section such that $s|_U \ne 0$
for some open set $U \subset Q$. Then for every
$s \in \Gamma^\infty_{\rm cpt} (L)$ with $\supp s \subset U$
there is a symbol $a \in \csym^{-\infty} (Q)$ such that
$\Op_0 (a) s_0 = s$.
\end {lemma}
\begin {proof}
Let $u \in \mathcal C^\infty_{\rm cpt} (Q)$ be the unique function with
$\supp u \subset U$ such that $us_0 = s$ and let $\sigma: TU \to \mathbb C$
be the unique smooth function satisfying
$s_0^\chi (v_x) = \sigma (v_x) s_0 (x)$. Then $\sigma|_{T_x U}$ has
compact support due to the choice of the cut-off function $\chi$
and fulfills $\sigma (0_x) = 1$ for all $x \in U$. Hence
\[
    S(x) := \frac{1}{(2\pi\hbar)^{n/2}} \int_{T_xQ}
            \left| \sigma (v_x) \right|^2 \mu_x(v_x)
    > 0
\]
is a strictly positive smooth function on $U$, so
$\widetilde b(v_x) := \frac{1}{S(x)} \cc{\sigma(v_x)}$ is a well-defined
smooth function on $TU$ such that $\widetilde b|_{T_xQ}$ has compact support
for every $x \in U$. Its fiberwise Fourier transform
\[
    b (\zeta_x) := \frac{1}{(2\pi\hbar)^{n/2}} \int_{T_xQ}
                   \eu^{\frac{\im}{\hbar} \brac{\zeta_x,v_x}}
                   \; \widetilde b (v_x) \; \mu_x(v_x)
\]
is an element of $\sym^{-\infty} (U)$. Finally
let $a = b (\pi^* u) \in \csym^{-\infty} (Q)$ and check by a
straightforward computation using the inverse fiberwise Fourier transform
that this $a$ fulfills the claim.
\end {proof}
\begin {lemma}
Let $s_0 \in \Gamma^\infty (L)$ be a section such that $s_0 (x_0) = 0$
and $s_0$ is transversal at $x_0$ for some $x_0 \in Q$. Then there
exists an open neighborhood $U$ of $x_0$ such that for every
$s \in \Gamma^\infty_{\rm cpt} (L)$ with $\supp s \subset U$
there exists a symbol $a \in \csym^{-\infty} (Q)$ satisfying
$\Op_0 (a)s_0 = s$.
\end {lemma}
\begin {proof}
First choose a contractible chart $V$ around $x_0$ with
coordinates $q^1, \ldots, q^n$ and induced coordinates
$v^1, \ldots, v^n$ resp.~$p_1, \ldots, p_n$ for the tangent resp.
cotangent vectors. Moreover, let $e \in \Gamma^\infty (L|_{V})$ be
a local non-vanishing section and let $\sigma_k : TV \to \mathbb C$
be the unique smooth functions with
$(\frac{\partial}{\partial v^k} s_0^\chi) (v_x) = \sigma_k (v_x) e(x)$
for $k = 1, \ldots, n$. Then $\sigma_k|_{T_x V}$ has compact support
for every $x \in V$. Due to the transversality of $s_0$ at $x_0$ we notice
that not all of the $\sigma_k (0_{x_0})$ vanish whence the function
\[
    S(x) := \frac{1}{(2\pi\hbar)^{n/2}} \sum_{k=1}^n
            \int_{T_xQ} \left| \sigma_k (v_x) \right|^2 \mu_x(v_x)
    \ge 0
\]
is smooth and non-negative over $V$ and satisfies $S (x_0) > 0$
due to the transversality of $s_0$ at $x_0$. Hence we can find an
open neighborhood $U\subset V$ of $x_0$ such that
over $U$ the function $S$ is strictly positive.
Now, let $s \in \Gamma^\infty_{\rm cpt} (L|_U)$ and
$u \in \mathcal C^\infty_{\rm cpt} (U)$ be the unique function
with $s = u e$. Furthermore define
$\widetilde b_k (v_x) = \frac{1}{S(x)} \cc{\sigma_k (v_x)}$ for
$v_x \in T_x U$ and $x \in U$. Then $\widetilde b_k$ is smooth
and $\widetilde b_k|_{T_x U}$ has compact support.
Thus its fiberwise Fourier transform $b_k$ is an element
of $\sym^{-\infty} (U)$. Next define the global symbol
$a \in \csym^{-\infty} (Q)$ by
\[
    a (\zeta_x) =
    \frac{\im}{\hbar} \, u(x) \sum_{k=1}^n b_k(\zeta_x) p_k (\zeta_x) .
\]
Then another straightforward computation using the inverse fiberwise
Fourier transform and involving one partial integration shows that
$\Op_0 (a)s_0 = s$.
\end {proof}

\begin {proposition} \label {SymbolSurjProp}
Let $s_0 \in \Gamma^\infty (L)$ be a transversal section and
$s \in \Gamma^\infty_{\rm cpt} (L)$ arbitrary. Then there exists a
symbol $a \in \csym^{-\infty} (Q)$ such that $\Op_0 (a) s_0 = s$.
\end {proposition}
\begin {proof}
The claim follows by the above two lemmas and a partition of unity
argument analogously to the proof of Lemma~\ref {TransLem}.
\end {proof}

Now let us come to the GNS construction for pseudodifferential
operators. Let $\mathcal A \subset \End (\Gamma^\infty_{\rm cpt} (L))$
be the algebra generated (not topologically but only algebraically)
by pseudodifferential operators of the form $\Opw (a)$, where
$a \in \csym^{-\infty} (Q)$.
Then $\mathcal A$ is a $^*$-algebra with the $^*$-structure given by the
(formal) operator adjoint on the pre-Hilbert space $\Gamma^\infty_{\rm cpt} (L)$.
The $^*$-involution now fulfills $\Opw (a)^* = \Opw (\cc a)$ which
follows from the definition of $\Opw$ in Theorem \ref {ThmPseudoOpDef}.
Note that in general $\mathcal A$ has no unit
element. The aim is now to define a positive functional on this
algebra such that the GNS pre-Hilbert space (over $\mathbb C$) is isometric
to $\Gamma^\infty_{\rm cpt} (L)$ and such that the induced GNS representation
coincides with the usual action of $\mathcal A$ on the sections of $L$.
So choose a transversal section $s_0 \in \Gamma^\infty (L)$ and set
\begin{equation} \label{TransNonformalOmegaDef}
    \omega_{s_0} (A) = \brac{s_0, As_0}_{L, \mu}, \qquad A \in \mathcal A.
\end{equation}
Then $\omega_{s_0}: \mathcal A \to \mathbb C$ comprises a well-defined
functional due to the fact that every operator
$\Opw (a)$ maps $s_0$ onto some section with compact support.
Then the following lemma is obvious.
\begin {lemma}
The so-called expectation value functional $\omega_{s_0}$ is
positive with Gel'fand ideal given by
$\mathcal J_{s_0} = \{ A \in \mathcal A \; | \; As_0 = 0\}$.
\end{lemma}
\begin{remark} In case $\mathcal A$ were the algebra of
bounded operators on a Hilbert space and $s_0$ a non-zero vector,
the analogue of the above functional would reproduce the Hilbert space
as GNS representation space and the
usual action of operators as GNS representation. In our case the result
is not quite as easy to achieve. The main problem lies in the proof of
the surjectivity of the isomorphism in the following theorem.
\end{remark}
\begin{theorem}
Let $s_0 \in \Gamma^\infty (L)$ be a transversal section and
$\omega_{s_0} : \mathcal A \to \mathbb C$ its expectation value functional.
Then the quotient pre-Hilbert space
$\mathfrak H_{s_0} = \mathcal A \big/ \mathcal J_{s_0}$ is isometric
to $\Gamma^\infty_{\rm cpt} (L)$ by the unitary map
\begin {equation} \label {SymbolUnitary}
    \mathfrak H_{s_0} \ni
    \psi_A \mapsto As_0
    \in \Gamma^\infty_{\rm cpt} (L),
    \qquad A \in \mathcal A.
\end {equation}
The GNS representation of $\mathcal A$
on $\mathcal A \big/ \mathcal J_{s_0}$ is unitarily equivalent to
the action of $\mathcal A$ on $\Gamma^\infty_{\rm cpt} (L)$ under
the unitary map (\ref {SymbolUnitary}).
\end {theorem}
\begin {proof}
First note that (\ref {SymbolUnitary}) is well-defined indeed
and isometric, hence injective. The surjectivity of
(\ref {SymbolUnitary}) is less trivial but follows from
Prop.~\ref {SymbolSurjProp} and relation (\ref {OpOpoCCNopo})
for $\order = 1/2$. Then the unitary equivalence of the GNS
representation with the  usual action of $\mathcal A$ on
$\Gamma^\infty_{\rm cpt} (L)$ follows immediately.
\end {proof}

\appendix
\section{Positive functionals and formal GNS construction}
\label{GNSApp}
For the reader's convenience let us recall some basic facts on positive
functionals and GNS constructions in deformation quantization. For a
detailed exposition and proofs we refer the reader to \cite{BW98a} and
\cite[App.~A]{BNW97b}.

The main observation which leads to the notion of a formal
GNS representation is the fact that $\mathbb R[[\lambda]]$ is an
\emph{ordered ring} by the following simple definition: a nonzero element
$a = \sum_{r=0}^\infty \lambda^r a_r \in \mathbb R[[\lambda]]$ is called
\emph{positive} if $a_{r_0} > 0$, where
$r_0 := \min\{r\in \mathbb N \; | \; a_r \ne 0 \}$. As easily verified
the sum and product of two positive elements is again positive and
each element is either negative, null, or positive. Moreover, after extending
complex conjugation coefficientwise to the ring
$\mathbb C[[\lambda]] = \mathbb R[[\lambda]]\oplus\im\mathbb R[[\lambda]]$
the following relations and inequalities hold: for
$z \in \mathbb C[[\lambda]]$ one has $\cc z z \ge 0$, and $\cc z z = 0$ is true
if and only if $z=0$. Note that we regard $\lambda$ as real, i.e.~we set
$\cc \lambda = \lambda$.

Now let $M$ be a symplectic or Poisson manifold with star product $*$ such
that the pointwise complex conjugation $f \mapsto \cc f$ is an
anti-automorphism of $*$, i.e.~we have $\cc{f *g} = \cc g * \cc f$ for
all $f, g \in \mathcal C^\infty (M)[[\lambda]]$. Then a
$\mathbb C[[\lambda]]$-linear functional
$\omega: \mathcal C^\infty(M)[[\lambda]] \to \mathbb C[[\lambda]]$ is
called \emph{positive} if and only if $\omega (\cc f * f) \ge 0$. Note that one
considers $\mathbb C[[\lambda]]$-linear and $\mathbb C[[\lambda]]$-valued
functionals instead of $\mathbb C$-linear and $\mathbb C$-valued ones. The
crucial point is that the notion of positivity extends naturally to this
framework. One immediately verifies now the Cauchy-Schwarz inequality for
positive functionals
\begin {equation} \label {CSU}
    \begin {array} {c}
        \omega (\cc f * g) = \cc{\omega (\cc g * f)} \\
        \omega (\cc f * g) \cc {\omega (\cc f * g)}
        \le
        \omega (\cc f * f) \omega (\cc g * g),
    \end {array}
\end {equation}
where $f, g \in \mathcal C^\infty (M)[[\lambda]]$. This implies by pure
algebraic reasoning as in the well-known case of $C^*$-algebras that
\begin {equation} \label {GelfandIdealDef}
    \mathcal J_\omega
    := \left\{ f \in \mathbb C^\infty (M)[[\lambda]] \; \big| \;
     \omega (\cc f * f) = 0 \right\}
\end {equation}
is a left-ideal, the so-called \emph{Gel'fand ideal} of $\omega$ in
$\mathcal C^\infty (M)[[\lambda]]$.

Given a positive $\mathbb C[[\lambda]]$-linear functional $\omega$ with
Gel'fand ideal $\mathcal J_\omega$ one considers the quotient
$\mathfrak H_\omega := C^\infty (M)[[\lambda]] \big/ \mathcal J_\omega$
which carries both a $\mathcal C^\infty (M)[[\lambda]]$-left module
structure and a Hermitian product: one defines the Hermitian product
by
\begin {equation} \label {GNSProductDef}
    \brac{\psi_f, \psi_g} := \omega (\cc f * g) \in \mathbb C[[\lambda]],
\end {equation}
where $\psi_f, \psi_g \in \mathfrak H_\omega$ denote the equivalence
classes of $f, g$. Indeed $\brac{\cdot, \cdot}$ is
$\mathbb C[[\lambda]]$-linear in the second argument and satisfies
$\brac{\psi_f, \psi_g} = \cc{\brac{\psi_g, \psi_f}}$ as well as the
positivity requirement $\brac{\psi_f, \psi_f} \ge 0$ for all
$\psi_f \in \mathfrak H_\omega$ and the non-degeneracy
$\brac{\psi_f, \psi_f} = 0$ if and only if $\psi_f = 0$. Thus $\mathfrak H_\omega$
becomes with this $\mathbb C[[\lambda]]$-valued Hermitian product a
pre-Hilbert space over the ring $\mathbb C[[\lambda]]$. Moreover the
$\mathcal C^\infty (M)[[\lambda]]$-left module structure induces the
\emph{GNS representation} $\pi_\omega$ of
$\mathcal C^\infty(M)[[\lambda]]$ on $\mathfrak H_\omega$ given by
\begin {equation} \label {FormalGNSRepDef}
    \pi_\omega (f) \psi_g := \psi_{f*g},
\end {equation}
which turns out to be a representation of the star product algebra. This
representation is even a $^*$-representation, i.e.~we have
\begin {equation} \label {GNSRepstarRep}
    \brac{\psi_g, \pi_\omega (f) \psi_h} =
    \brac{\pi_\omega (\cc f)\psi_g, \psi_h}
\end {equation}
for all $\psi_g, \psi_h \in \mathfrak H_\omega$ and
$f \in \mathcal C^\infty (M)[[\lambda]]$. In this sense we also write
$\pi (f)^* = \pi (\cc f)$.

Finally we would like to mention that the domain of $\omega$ need not be the
whole algebra $\mathcal C^\infty (M)[[\lambda]]$ but rather a two-sided
ideal $\mathcal B$ stable under complex conjugation
$\cc {\mathcal B} = \mathcal B$. Then $\mathcal J_\omega$ turns out to be
even a left-ideal in $\mathcal C^\infty (M)[[\lambda]]$ due to the
Cauchy-Schwarz inequality. Thus the GNS representation $\pi_\omega$ of
$\mathcal B$ on $\mathfrak H_\omega = \mathcal B \big/ \mathcal J_\omega$
extends naturally to a $^*$-representation of the whole algebra
$\mathcal C^\infty (M)[[\lambda]]$, a situation which is present
throughout this paper. Here the two-sided ideal is given e.g.~by
$\mathcal C^\infty_Q (T^*Q)[[\lambda]]$.

For further generalizations to arbitrary ordered rings as well as more
results in the context of deformation quantization including the field of
formal Laurent series $\LS{\mathbb C}$ and its field extension
$\CNP{\mathbb C}$ we refer to $\cite{BW98a}$. For a detailed treatment
of the implementation of symmetries as unitary operators in GNS
representations, the `pull-back' of positive functionals and the
corresponding GNS representations we refer to \cite{BNW97b}.

%\bibliographystyle{wde}
%\bibliography{articles,books,preprints,ownpreps,misc}

\end {document}